\documentclass[10pt,a4paper]{article}
\usepackage[utf8]{inputenc}
\usepackage[english]{babel}
\usepackage[T1]{fontenc}
\usepackage{amsmath}
\usepackage{amsfonts}
\usepackage{bbm}
\usepackage{graphicx}
\usepackage{subcaption}
\usepackage{array}
\usepackage{bm}
\usepackage{authblk}
\usepackage{upgreek}
\usepackage[left=2.5cm,right=2.5cm,top=3cm,bottom=3cm]{geometry}
\usepackage{fancyhdr}
{ \fancyhead[L]{\nouppercase{\rightmark}}
\fancyhead[R]{\nouppercase{\leftmark}}  
 \fancyhead[R]{\textbf{ \rightmark}}
 \fancyhead[L]{\textbf{\thepage}}
 \fancyfoot[C]{}}

\usepackage{lipsum}           
\usepackage{adjustbox}
\usepackage{transparent}

\usepackage{isomath}
\usepackage{mathtools}
\usepackage{mathrsfs}
\usepackage{amsthm}
\usepackage{slashed}
\usepackage{amssymb}
\usepackage{enumitem}
\usepackage{hyperref}
\hypersetup{colorlinks=true, 
citecolor=red}
\usepackage[new]{old-arrows}

\theoremstyle{definition}
\newtheorem{mdef}{{Definition}}[section]

\theoremstyle{definition}
\newtheorem{mex}{Example}[section]

\theoremstyle{definition}
\newtheorem{mrmk}{{Remark}}[section]

\theoremstyle{plain}
\newtheorem{mth}{Theorem}[section]

\theoremstyle{plain}
\newtheorem{mlem}{{Lemma}}[section]

\theoremstyle{plain}
\newtheorem{mprop}{{Proposition}}[section]

\theoremstyle{plain}
\newtheorem{mcor}{{Corollary}}[section]

\theoremstyle{definition}
\newtheorem{mnot}{{Notation}}[section]

\newlength{\bibitemsep}
\newlength{\bibparskip}
\let\oldthebibliography\thebibliography
\renewcommand\thebibliography[1]{%
  \oldthebibliography{#1}%
  \setlength{\parskip}{\bibitemsep}%
  \setlength{\itemsep}{\bibparskip}%
}

\usepackage{tikz-cd} 
\usetikzlibrary{positioning, shapes, arrows.meta,decorations.markings,decorations.pathmorphing,shapes}
\tikzset{
    answer/.style={rectangle, draw, text width=15em, text badly centered, node distance=1cm, inner sep=0pt, minimum height=4em},
    block/.style={rectangle, draw, text width=10em, text centered},
     block2/.style={rectangle, draw, text width=5em, text centered},
     block3/.style={rectangle, draw, text width=5em, text centered, color=white},
}

\DeclareMathSymbol{\upLambda}{\mathalpha}{operators}{3}

\begin{document}
\title{Noncommutative geometry on the Berkovich projective line}
\author[1]{Masoud Khalkhali}
\author[2]{Damien Tageddine}
\affil[1]{Department of mathematics, The University of Western Ontario}
\affil[2]{Department of mathematics and statistics, McGill University}
\date{}

\maketitle

\begin{abstract}
\noindent
In this paper, we construct several $C^*$-algebras associated to the Berkovich projective line  $\mathbb{P}^1_{\mathrm{Berk}}({\mathbb{C}_p})$. In the commutative setting, we construct a spectral triple as a direct limit over finite $\mathbb{R}$-trees. More general $C^*$-algebras generated by partial isometries are also presented. We use their representations to associate a Perron-Frobenius operator and a family of projection-valued measures. Finally, we show that invariant measures, such as the Patterson-Sullivan measure, can be obtained as KMS-states of the crossed product algebra with a Schottky subgroup of $\mathrm{PGL}_2(\mathbb{C}_p)$.
\end{abstract}
{
  \hypersetup{linkcolor=black}
  \tableofcontents
}

\newpage 

\section{Introduction}
\noindent
The aim of this paper is to show that geometries over non-Archimedean fields provide natural examples of noncommutative geometries; we can refer to the work \cite{consani_spectral_2003,consani_non-commutative_2003} for similar approaches. We focus on the Berkovich projective line, \(\mathbb{P}^1_{\text{Berk}}(\mathbb{C}_p)\), as a fundamental example of Berkovich's theory. For further details on noncommutative geometry, we refer to the books \cite{connes_noncommutative_1994,connes_noncommutative_2007}; for an extensive exploration of the interplay between number theory and noncommutative geometry, we refer to the series of work \cite{connes_scaling_2015, connes_geometry_2016,connes_noncommutative_2007}.\\

\noindent
The primary objective of this work is to construct and analyze several C*-algebras and spectral triples that can be associated to the Berkovich projective line. These C*-algebras are expected to encapsulate some essential geometric and arithmetic features of \(\mathbb{P}^1_{\text{Berk}}(\mathbb{C}_p)\) in a manner analogous to how C*-algebras encode information about classical spaces.\\

\noindent
Let $K$ be an algebraically closed field that is complete with respect to a nontrivial non-Archimedean absolute value $|\cdot|$. 
In Section \ref{Sect1}, we begin by reviewing the definition of the Berkovich projective line $\mathbb{P}^1_{\mathrm{Berk}}(K)$ and its fundamental properties. In particular, after restriction to the case $K=\mathbb{C}_p$, we recall the classification of the points of $\mathbb{P}^1_{\mathrm{Berk}}(\mathbb{C}_p)$ into four types according to Berkovich's classification theorem. We focus on the $\mathbb{R}$-tree structure obtained from a projective limit of finite trees and refer to \cite{baker_potential_nodate} for more details. We will also introduce a hyperbolic metric structure and the definition of the hyperbolic space denoted by $\mathbb{H}_{\mathrm{Berk}}(\mathbb{C}_p)$.\\

\noindent
In Section \ref{Sect2}, we provide a first construction of a commutative spectral triple represented by the datum $(C_{\mathrm{Lip}}(\mathbb{P}^1_{\mathrm{Berk}}(\mathbb{C}_p)),\mathcal{H},D)$ as an inverse limit of finite spectral triples associated with finite trees. The result is summarized in Theorem \ref{theorem1}. We closely follow the general construction of an inverse limit of spectral triples given in \cite{floricel_inductive_2017}. \\

\noindent
In Section \ref{Sect4}, we propose an alternative construction of $C^*$-algebra associated to $\mathbb{P}^1_{\mathrm{Berk}}(\mathbb{C}_p)$ relying on the identification of the projective line with the Wa\.zewski universal dendrite introduced in Section \ref{Sect3}. The $C^*$-algebra $\mathcal{O}_{\mathbb{P}^1_{\mathrm{Berk}}(\mathbb{C}_p)}$ is generated by partial isometries indexed by the branching points of the universal dendrite. It is in fact the full shift $C^*$-algebra associated to a countable alphabet, with letters given by the rational numbers in $(0,1)$. This type of algebra bears similar properties than the Cuntz-Krieger algebras, see for instance \cite{marcolli_cuntzkrieger_2011, jorgensen_states_2011}. The results are summarized in Proposition \ref{prop1} and Theorem \ref{mth2}. Using the representation space, we also define a Perron-Frobenius operator and projection-valued measures.\\

\noindent
Finally, in Section \ref{Sect5}, we use the fact that $\mathrm{PGL}_2(\mathbb{C}_p)$ is the isometry group of $\mathbb{P}^1_{\mathrm{Berk}}(\mathbb{C}_p)$. We exhibit a unitary representation of $\mathrm{PGL}_2(\mathbb{C}_p)$ through the left action on the projective line. This allows us to identify the boundary of $\mathbb{P}^1_{\mathrm{Berk}}(\mathbb{C}_p)$  as the limit set of the action. The boundary coincides with the $p$-adic .We construct the crossed product $C^*$-algebras, $C_{\mathrm{Lip}}(\mathbb{P}^1(\mathbb{C}_p))\rtimes \mathrm{\Gamma}$ with a Schottky group. The Patterson-Sullivan measure is obtained as a KMS-state in Theorem \ref{thm3}.\\

\noindent
\noindent
The long-term goal of the project we start in this paper is to provide a unified approach to the noncommutative geometry treatment of number theory in both the Archimedean and non-Archimedean settings, as studied in the papers \cite{connes_scaling_2015,connes_geometry_2016,consani_spectral_2003,consani_non-commutative_2003}. This follows the general philosophy of Berkovich aiming at describing analogues of complex analytic spaces in the complex $p$-adic case, providing a $p$-adic analytic space where such a theory can be defined.
Methods of noncommutative geometry, specially spectral triples and Bost-Connes systems, have proved to be powerful to capture the geometry of such \textit{wild spaces} that often appear in number theory \cite{connes_physics_nodate}.\\
Apart from potential applications to number theory, the class of $C^*$-algebras introduced in this paper and that appears in the potential theory and dynamics of rational maps on the Berkovich projective line provide new examples of infinite graph $C^*$-algebras. In a different direction, we would like to mention that the present work will hopefully provide a framework to test and rigorously define a  version of the AdS/CFT correspondence along the lines of \cite{heydeman_tensor_2017,gubser_edge_2017}.
\section{The Berkovich projective line}
\label{Sect1}
In this section, we recall the definition of \textit{multiplicative seminorms} and introduce the construction of Berkovich spaces.
\subsection{The multiplicative seminorms and Berkovich spaces}
\begin{mdef}
Let $K$ be a field with absolute value $|\cdot|$ and let $A$ be a $K$-algebra. A \textit{multiplicative seminorm} on $A$ is a map
\begin{equation*}
A\rightarrow \mathbb{R}_{\geq 0}, \quad a\mapsto \|a\|
\end{equation*}
such that
\begin{itemize}[leftmargin=0.6cm]
\setlength\itemsep{0em}
\item[(1)] $\|\cdot\|$ restricts to $|\cdot|$ on $K$;
\item[(2)] $\|a+b\|\leq \|a\| + \|b\|$ for all $a,b\in A$;
\item[(3)] $\|ab\|=\|a\|\cdot \|b\|$ for all $a,b\in A$.
\end{itemize}
If in addition $a=0$ is the only element with $\|a\|=0$, then $\|\cdot\|$ is a \textit{multiplicative norm} (which is the same as an absolute value on $A$ extending $|\cdot|$). We call $\mathrm{ker}\|\cdot \|=\lbrace a\in A : \|a\|=0\rbrace$ the kernel of $\|\cdot\|$; it is a prime ideal of $A$.\\ \noindent
A $K$-algebra $A$ with a fixed multiplicative norm such that $A$ is complete with respect to this norm is a \textit{Banach algebra} over $K$.
\end{mdef}
\begin{mdef}
Let $(K,|\cdot|)$ be a field with absolute, let $a\in K$ and $r\geq 0$.
Then, 
\begin{equation*}
D(a,r):=\lbrace \xi \in K : |\xi-a|\leq r\rbrace
\end{equation*}
is the closed disk of radius $r$ around $a$.
\end{mdef}
\begin{mex}
Let $(K,|\cdot|)$ be a complete non-archimedean field.
\begin{itemize}[leftmargin=0.6cm]
\setlength\itemsep{0em}
\item[(1)] For any $a\in K$, the map
\begin{equation*}
f\mapsto \|f\|_{a,0}:=|f(a)|
\end{equation*}
is a multiplicative seminorm on $K[x]$.
\item[(2)] For any $a\in K$ and any $r>0$, the map
\begin{equation*}
f\mapsto \|f\|_{a,r}:=|f(x+a)|_r
\end{equation*}
is a multiplicative (semi)norm on $K[x]$.
\item[(3)] Let $(a_n)$ be a sequence in $K$ and $(r_n)$ a strictly decreasing sequence in $\mathbb{R}_{> 0}$ such that $D(a_{n+1},r_{n+1})\subset D(a_n,r_n)$ for all $n$. Then 
\begin{equation*}
f\mapsto \|f\|=\lim_{n\rightarrow \infty}\|f\|_{a_n,r_n}
\end{equation*}
is a multiplicative seminorm on $K[x]$.
\end{itemize}
\end{mex}
\begin{mdef}
Let $K$ be a complete and algebraically closed non-Archimedean field and let $A$ be a finitely generated $K$-algebra, so that $A$ is the coordinate ring of the affine $K$-variety $X=\mathrm{Spec}(A)$. Then the \textit{Berkovich space} associated to $A$ or $X$ is
\begin{equation}
\mathrm{Berk}(A):=X_\mathrm{Berk}:=\lbrace\|\cdot\|:A\rightarrow \mathbb{R} \ \text{multiplicative seminorm on} A\rbrace
\end{equation}
the set of multiplicative seminorms on $A$. The topology on $X_\mathrm{Berk}$ is the Gelfand topology, i.e. the weakest topology that makes the maps $X_\mathrm{Berk}\rightarrow \mathbb{R}$, $\|\cdot\|\rightarrow \|f\|$, continuous for all $f\in A$. (Concretely, this means that any open set is a union of finite intersections of sets of the form $U_{f,a,b}=\lbrace\|\cdot\|\in X_\mathrm{Berk} : a<\|f\|<b\rbrace$.)
\end{mdef}
\noindent
Elements of $X_\mathrm{Berk}$ are called points. For a point $\xi\in X_\mathrm{Berk}$, the corresponding seminorm is denoted $\|\cdot\|_{\xi}$.
\subsection{Berkovich's classification theorem}
\begin{mth}[Berkovich's Classification Theorem] Every point $x\in \mathbb{A}^1_{\mathrm{Berk}}(K)$ corresponds to a nested sequence $D(a_1,r_1)\supseteq D(a_2,r_2)\supseteq D(a_3,r_3)\supseteq \cdots$ of closed disks, in the sense
\begin{equation*}
\|f\|_x=\lim_{n\rightarrow \infty}\|f\|_{D(a_n,r_n)}
\end{equation*}
Two such nested sequences define the same point of $\mathbb{A}^1_{\mathrm{Berk}}(K)$ if and only if
\begin{itemize}
\setlength\itemsep{0em}
\item[(a)] each has a nonempty intersection, and their intersections are the same; or
\item[(b)] both have empty intersection, and the sequences are cofinal.
\end{itemize}
\end{mth}
\noindent
This brings us to Berkovich's classification of elements of $\mathbb{A}^1_{\mathrm{Berk}}(K)$ into four types of points according to the nature of $D(a,r)=\cap \ D(a_n,r_n)$:
\begin{itemize}
\setlength\itemsep{0em}
\item[(I)] $D(a,r)$ is a point ($r=\lim r_i=0$).
\item[(II)] $D(a,r)$ is a closed disk of radius $r=\lim r_i>0$ belonging to the value group $|K^*|$ of $K$. The corresponding discs $D(a,r)$ are called \textit{'rational'}.
\item[(III)] $D(a,r)$ is a closed disk of radius $r=\lim r_i>0$ that does not belong to the value group $|K^*|$ of $K$. The corresponding discs $D(a,r)$ are called \textit{'irrational'}.
\item[(IV)] $D(a,r)$ is empty. As noted above, necessarily $\lim r_i >0$.
\end{itemize}
\begin{mrmk}
Berkovich's classification theorem is crucial to understand the action of the automorphism group of the tree. Specifically, the type of a point is invariant under the action of the group. Moreover, type II points are dense in the tree, therefore one can understand the action by restriction to the branching points.
\end{mrmk}
\subsection{The Berkovich affine line $\mathbb{A}^1_{\mathrm{Berk}}(\mathbb{C}_p)$}
In this section and in the rest of the paper, we consider (unless stated otherwise) the Berkovich affine line over $\mathbb{C}_p$, $\mathbb{A}^1_{\mathrm{Berk}}(\mathbb{C}_p)=\mathrm{Berk}\left( \mathbb{C}_p\left[ x\right] \right) $. According to the Berkovich classification of points, the affine line $\mathbb{A}^1_{\mathrm{Berk}}(\mathbb{C}_p)$ possesses four types of points. Type I points are the classical points $\mathbb{A}^1(\mathbb{C}_p)=\mathbb{C}_p$. Type II and III correspond to closed disks $D(a,r)$ with $r>0$ in, respectively, not in, the value group $p^\mathbb{Q}$ of $\mathbb{C}_p$. We will use the notation $\zeta_{a,r}$ for these points. By extension, we will identify the points $a$ of type I with $\zeta_{a,0}$.\\

\noindent
The Berkovich affine line is endowed with a natural partial order. For any two points $\xi,\xi'\in \mathbb{A}^1_{\mathrm{Berk}}(\mathbb{C}_p)$, we say that:
\begin{equation*}
\xi\leq \xi' \quad \text{if and only if} \quad \forall f\in \mathbb{C}_p[x] : \ \|f\|_{\xi}\leq \|f\|_{\xi'}
\end{equation*}
In particular, for two disks $D_1$ and $D_2$, $\|\cdot\|_{D_1}\leq \|\cdot\|_{D_2}$ if and only if $D_1\subset D_2$. The relation $\leq$ defines a partial order on the affine line. For each pair of points $\xi,\xi'\in \mathbb{A}^1_{\mathrm{Berk}}(\mathbb{C}_p)$, there is a unique least upper bound $\xi\vee \xi'\in \mathbb{A}^1_{\mathrm{Berk}}(\mathbb{C}_p)$ with respect to this partial order.\\
 Moreover, if $\xi\leq \xi'$, then we write $\left[ \xi,\xi'\right] $ for the set of points $\xi''$ such that $\xi\leq \xi''\leq \xi'$. If the points $\xi$ and $\xi'$ are of type I, II or III, then the set $\left[ \xi,\xi'\right] $ is homeomorphic to a closed interval in $\mathbb{R}$.
\begin{mprop}
The space $\mathbb{A}^1_{\mathrm{Berk}}(\mathbb{C}_p)$ is Hausdorff, locally compact and uniquely path-connected.
\end{mprop}
\noindent
We can now introduce different metrics on the Berkovich space $\mathbb{A}^1_{\mathrm{Berk}}(\mathbb{C}_p)$.
\begin{mdef}
The diameter of a point $\xi\in \mathbb{A}^1_{\mathrm{Berk}}(\mathbb{C}_p)$ is 
\begin{equation*}
\mathrm{diam}(\xi):=\inf \left\lbrace \|x-a\|_\xi : a\in \mathbb{C}_p \right\rbrace 
\end{equation*}
\end{mdef}
\noindent
Therefore, points of types I, II and III have a diameter $\mathrm{diam}(\zeta_{a,r})=r$, and for a type IV point $\xi$ represented by a nested sequence of disks $D(a_n,r_n)$, we have $\mathrm{diam}(\xi)=\lim_{n\rightarrow \infty}r_n>0$.\\

\noindent
There are two metrics that one can define using the $\mathrm{diam}$ map.
\begin{mdef}
For any couple of points $\xi,\xi'\in \mathbb{A}^1_{\mathrm{Berk}}(\mathbb{C}_p)$, we define the \textit{small metric} by 
\begin{equation*}
d(\xi,\xi')=2 \ \text{diam}(\xi\vee \xi')-\mathrm{diam}(\xi)-\mathrm{diam}(\xi')
\end{equation*}
In the case that $\xi,\xi'$ are both not of type I, the \textit{big metric} is defined by
\begin{equation*}
\rho(\xi,\xi')=2\log_v\mathrm{diam}(\xi\vee \xi') - \log_v \mathrm{diam}(\xi) - \log_v \mathrm{diam} (\xi')
\end{equation*}
\end{mdef}
\noindent
The small metric measures the total variation of diameter as we move along the path from $\xi$ to $\xi'$. The big metric is the analogue of the hyperbolic metric on the upper half plane $ds^2=(dx^2+dy^2)/y^2$ in $\mathbb{C}$. In fact, the \textit{hyperbolic part} of $\mathbb{A}^1_{\mathrm{Berk}}(\mathbb{C}_p)$ consists of points of types II, III and IV. The big metric is invariant under the action of $\mathrm{PGL}_2(\mathbb{C}_p)$, analogous to the fact that the standard hyperbolic metric is invariant under the action of $\mathrm{PSL}_2(\mathbb{R})$. We can now introduce the tree-like structure of the Berkovich line.
\begin{mdef}[$\mathbb{R}$-tree]
An $\mathbb{R}$-tree is a metric space $(T,d)$ such that for any two points $x,y\in T$, there exists a unique path $\left[ x,y\right] $ in $T$ joining $x$ to $y$, which is a geodesic segment; i.e. the map $\gamma : \left[ a,b\right] \rightarrow T$ giving the path can be chosen so that $d(\gamma(u),\gamma(v))=|u-v|$ for all $u,v\in [a,b]$.\\

\noindent
A point $x\in T$ is a \textit{branch point} if $T\backslash \left\lbrace x \right\rbrace$ has at least three connected components in the metric topology. A point $x$ is an \textit{endpoint} if $T\backslash \left\lbrace x \right\rbrace$ is connected. An \textit{ordinary} point $x$ is a point for which $T\backslash \left\lbrace x \right\rbrace$ has exactly two connected components.\\

\noindent
The metric topology on $T$ is called the \textit{strong topology}. The \textit{weak topology} is defined using the \textit{tangent direction} at $x\in T$. A \textit{tangent direction} is an equivalence class of path $[x,y]$ with $x\neq y$, where two paths are equivalent if they share an initial segment. The connected components of $T\backslash \left\lbrace x \right\rbrace$ are in one-to-one correspondence with the tangent directions at $x$. For a tangent direction $v$ at $x$, we define $U(x,v)=\lbrace y\in T : y\neq x, [x,y]\in v\rbrace$. The sets $U(x,v)$ is a subbasis for the \textit{weak topology}.
\end{mdef}
\begin{mth}The following statements hold:
\begin{itemize}[leftmargin=0.6cm]
\vspace{-1em}
\setlength\itemsep{0em}
\item[~]
\item[(1)]The Berkovich line $(\mathbb{A}^1_{\mathrm{Berk}}(\mathbb{C}_p),d)$ equipped with the small metric $d$ has a $\mathbb{R}$-tree structure. The Berkovich topology on $\mathbb{A}^1_{\mathrm{Berk}}(\mathbb{C}_p)$ is the weak topology on this $\mathbb{R}$-tree.
\item[(2)]The space $(\mathbb{A}^1_{\mathrm{Berk}}(\mathbb{C}_p)\backslash\mathbb{C}_p,\rho)$ equipped with the big metric $\rho$ has a $\mathbb{R}$-tree structure. The subspace topology on $\mathbb{A}^1_{\mathrm{Berk}}(\mathbb{C}_p)\backslash\mathbb{C}_p$ is the weak topology on this $\mathbb{R}$-tree.
\end{itemize}
\end{mth}
\begin{figure}
\begin{subfigure}{.5\textwidth}
  \centering
  \adjustbox{trim=3cm 9cm 3cm 0cm}{
\resizebox{120mm}{!}{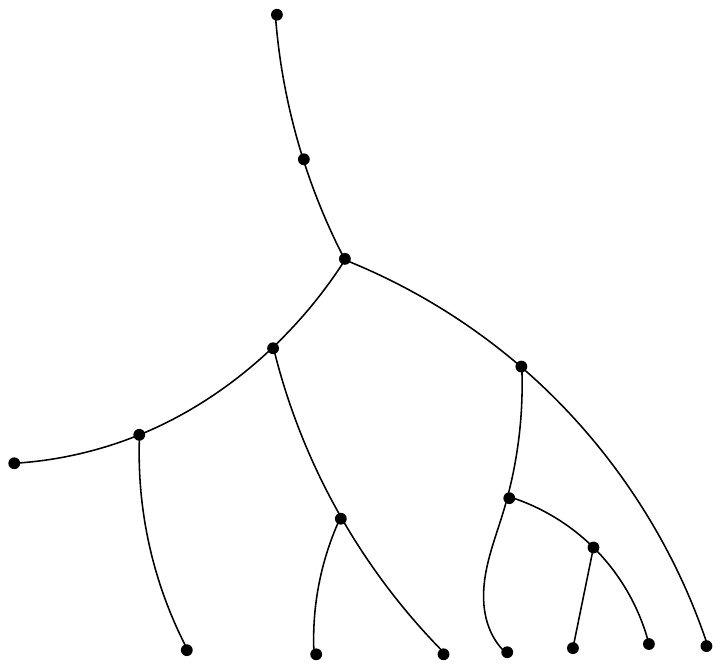}}
\end{subfigure}%
\begin{subfigure}{.5\textwidth}
  \centering
  \includegraphics[width=.8\linewidth]{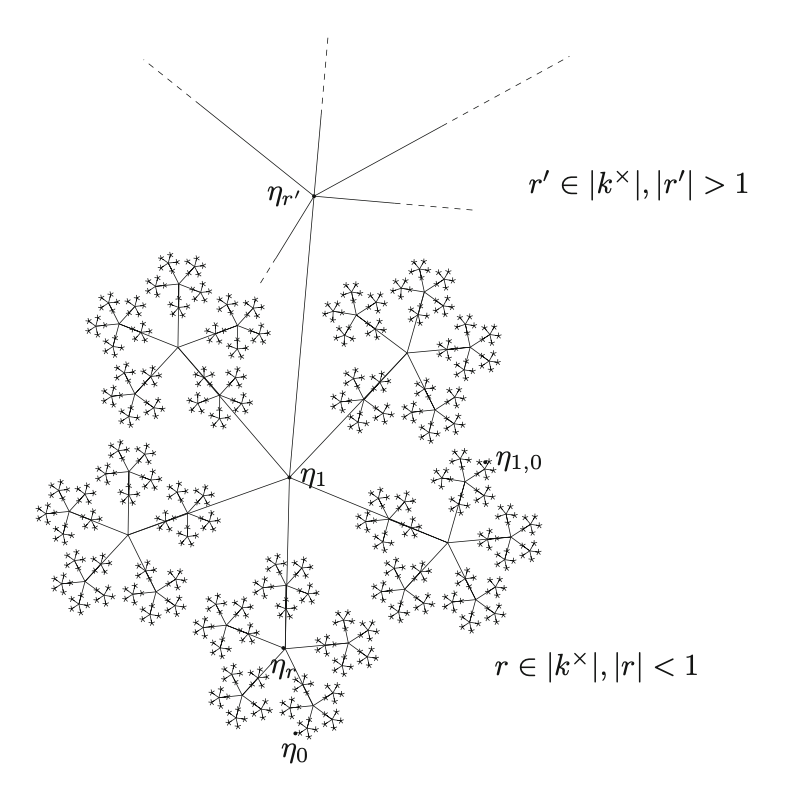}
\end{subfigure}
\caption{The Berkovich affine line. Picture on the right is taken from \cite{Poineau2021}}
\label{fig:fig}
\end{figure}
\subsection{The Berkovich projective line $\mathbb{P}^1_{\mathrm{Berk}}(\mathbb{C}_p)$ and the $\mathbb{R}$-tree structure}
As a set, the Berkovich projective line, denoted by $\mathbb{P}^1_{\mathrm{Berk}}(\mathbb{C}_p)$, is obtained by adding a type I point, denoted $\infty$ to the affine line $\mathbb{A}^1_{\mathrm{Berk}}(\mathbb{C}_p)$. The set $\mathbb{P}^1_{\mathrm{Berk}}(\mathbb{C}_p)$ is then equipped with the one-point compactification topology.\\
\noindent
The Berkovich hyperbolic line $\mathbb{H}_{\mathrm{Berk}}(\mathbb{C}_p)$ is a subset of $\mathbb{P}^1_{\mathrm{Berk}}(\mathbb{C}_p)$ consisting of all points of type II, III, and IV. One also defines $\mathbb{H}^\mathbb{Q}_{\mathrm{Berk}}(\mathbb{C}_p)$, the set of all type II points, and $\mathbb{H}^\mathbb{R}_{\mathrm{Berk}}(\mathbb{C}_p)$ for the set of points of type II and III.
\begin{mprop}
The subset $\mathbb{H}^\mathbb{Q}_{\mathrm{Berk}}(\mathbb{C}_p)$ is dense in $\mathbb{P}^1_{\mathrm{Berk}}(\mathbb{C}_p)$.
\end{mprop}
\noindent
The projective line $\mathbb{P}^1_{\mathrm{Berk}}(\mathbb{C}_p)$ also inherits a $\mathbb{R}$-tree structure. It is the compactification of $\mathbb{A}^1_{\mathrm{Berk}}(\mathbb{C}_p)$ seen as an $\mathbb{R}$-tree.
Following the description of Baker and Rumely \cite{baker_potential_nodate}, one can \textit{navigate} on the Berkovich projective line as follows. One starts from the so-called Gauss point $\zeta_{0,1}$, corresponding to the unit disk which is a type II point, and chooses between infinitely many (countable) branches in which to travel; there are one branch for each element of the residue field $\overline{\mathbb{F}}_p$. On the chosen direction, at each point of type II there are again infinitely many new branches in which to travel.\\

\noindent
We can classify the points in the $\mathbb{R}$-tree structure as follows:
\begin{itemize}[leftmargin=0.6cm]
\setlength\itemsep{0em}
\item[(1)] Points of types I and IV are endpoints.
\item[(2)] Points of type II are the only branching points, with infinite (countable) tangent directions in bijection with the residue field $\mathbb{P}^1(\overline{\mathbb{F}}_p)$.
\item[(3)] Points of type III are ordinary points.
\end{itemize}
\noindent 
There is an important description of the Berkovich projective line as a \textit{profinite $\mathbb{R}$-tree,} i.e. an inverse limit of $\mathbb{R}$-trees. Indeed, $\mathbb{P}^1_{\mathrm{Berk}}(\mathbb{C}_p)$ is homeomorphic to the inverse limit $\lim\limits_{\leftarrow}\mathrm{\Gamma}$ over all finite $\mathbb{R}$-trees $\mathrm{\Gamma}\subset \mathbb{P}^1_{\mathrm{Berk}}(\mathbb{C}_p)$.
More precisely, let $S=\lbrace D(a_1,r_1),\dots, D(a_n,r_n)\rbrace$ be a finite set of (rational or irrational) discs of positive radius contained in $D(0,1)$. For simplicity, we assume that $D(0,1)\notin S$. To each disc $D(a_i,r_i)$, there is an associated point $\zeta_{a_i,r_i}\in \mathbb{A}^1_{\mathrm{Berk}}(\mathbb{C}_p)$, which is a point of type II or III.\par 
\noindent
Define the \textit{graph of discs} $\mathrm{\Gamma}_S$ to be the union of the associated lines of discs $\left[ \zeta_{a_i,r_i},\zeta_{\mathrm{Gauss}}\right] $:
\begin{equation*}
\mathrm{\Gamma}_S=\bigcup_{i=1}^n\left[ \zeta_{a_i,r_i},\zeta_{\mathrm{Gauss}}\right]
\end{equation*}
If $S_1$ and $S_2$ are any two finite sets of discs, then $\Gamma_{S_1}$ and $\mathrm{\Gamma}_{S_2}$ are both $\mathbb{R}$-subtrees of $\mathrm{\Gamma}_{S_1\cup S_2}$. Moreover, the embedding of  $\mathrm{\Gamma}_{S_i}$ in $\mathrm{\Gamma}_{S_1\cup S_2}$ is an isometry (with respect to either metric). Let $\mathcal{F}$ be the collection of all finite graphs of the form $\mathrm{\Gamma}_S$ as above. Then, $\mathcal{F}$ is a directed set under inclusion and we write $\mathrm{\Gamma}\leq \mathrm{\Gamma}'$ if $\mathrm{\Gamma}\subseteq \mathrm{\Gamma}'$ as subsets of $D(0,1)$. Thus, whenever $\mathrm{\Gamma} \leq \mathrm{\Gamma}'$ there is an inclusion map $i_{\mathrm{\Gamma},\mathrm{\Gamma}'}:\mathrm{\Gamma}\rightarrow \mathrm{\Gamma}'$.\\

\noindent 
Similarly, write $\mathcal{F}^{\mathbb{Q}}$ (resp. $\mathcal{F}^{\mathbb{R}}$) for the subset of $\mathcal{F}$ consisting of graphs which are union of all arcs connecting two points in $\mathbb{H}^{\mathbb{Q}}_{\mathrm{Berk}}$ (resp. $\mathbb{H}^{\mathbb{R}}_{\mathrm{Berk}}$).\\

\noindent
There is also a retraction map $r_{\mathrm{\Gamma}',\mathrm{\Gamma}}:\mathrm{\Gamma}'\rightarrow \mathrm{\Gamma}$ defined whenever $\mathrm{\Gamma}\leq \mathrm{\Gamma}'$. This is a general property of $\mathbb{R}$-trees: since there is a unique path between any two points of $\mathrm{\Gamma}'$, if $x\in \Gamma'$ we can define $r_{\mathrm{\Gamma}',\mathrm{\Gamma}}(x)$ to be, for any $y\in \mathrm{\Gamma}$, the first point where the unique path in $\mathrm{\Gamma}'$ from $x$ to $y$ intersects $\mathrm{\Gamma}$. This definition is independent of the choice of $y$, and one sees from the definition that $r_{\mathrm{\Gamma}',\mathrm{\Gamma}}(x)=x$ if and only if $x\in \Gamma$. In particular, $r_{\mathrm{\Gamma}',\mathrm{\Gamma}}$ is surjective.
\begin{mth}
There is a canonical homeomorphism
\begin{equation*}
\mathbb{P}^1_{\mathrm{Berk}}(\mathbb{C}_p)\simeq \lim_{\xleftarrow[\mathrm{\Gamma}\in \mathcal{F}]{}}\mathrm{\Gamma}
\end{equation*}
Then $\mathbb{P}^1_{\mathrm{Berk}}(\mathbb{C}_p)$ is a compact Hausdorff space identified with the inverse limit of its finite subgraphs with the inverse limit topology.
\end{mth}
\noindent
We can also consider the direct limit with respect to the inclusion maps $i_{\mathrm{\Gamma},\mathrm{\Gamma}'}$. In this case, the limit is isomorphic (as a set) to $\mathbb{H}_{\mathrm{Berk}}(\mathbb{C}_p)$.
\begin{mth}
There are canonical bijections 
\begin{align*}
&\mathbb{H}_{\mathrm{Berk}}(\mathbb{C}_p)\simeq \lim_{\xrightarrow[\mathrm{\Gamma}\in \mathcal{F}]{}}\mathrm{\Gamma}\\
\mathbb{H}^{\mathbb{Q}}_{\mathrm{Berk}}(\mathbb{C}_p)\simeq &\lim_{\xrightarrow[\mathrm{\Gamma}\in \mathcal{F}^{\mathbb{Q}}]{}}\mathrm{\Gamma}, \quad \mathbb{H}^{\mathbb{R}}_{\mathrm{Berk}}\simeq \lim_{\xrightarrow[\mathrm{\Gamma}\in \mathcal{F}^{\mathbb{R}}]{}}\mathrm{\Gamma}
\end{align*}
\end{mth}
\noindent
In other words, the branching points $\mathbb{H}^{\mathbb{Q}}_{\mathrm{Berk}}(\mathbb{C}_p)$ is identified with the following set 
\begin{equation*}
\mathbb{H}^{\mathbb{Q}}_{\mathrm{Berk}}(\mathbb{C}_p)=\left\lbrace (\zeta_{a_1,r_1},\zeta_{a_2,r_2},\dots \zeta_{a_n,r_n}, \zeta_{a_n,r_n},\dots) : \text{s.t. $\zeta_{a_i,r_i}\in \mathrm{Type\ II}$ for all $i$} \right\rbrace.
\end{equation*}
We refer to \cite{baker_potential_nodate} for more details on this construction. 
\section{A commutative spectral triple on $\mathbb{P}^1_{\mathrm{Berk}}(\mathbb{C}_p)$}
\label{Sect2}
\subsection{Spectral triples on finite $\mathbb{R}$-trees}
Consider a finite $\mathbb{R}$-tree $\mathrm{\Gamma}\in \mathcal{F}$ such that $\mathrm{\Gamma}=(\mathcal{V},\mathcal{E})$. For any a vertex $v$, denote by $n_v$ the number of adjacent vertices to $v$. We construct a spectral triple over $\mathrm{\Gamma}$, denoted by $(A_\mathrm{\Gamma},\mathcal{H}_\mathrm{\Gamma},D_\mathrm{\Gamma})$. Let $A_\mathrm{\Gamma}$ be the *-algebra of (Lipschitz continuous functions) $C_{\mathrm{Lip}}(\mathrm{\Gamma})$ on $(\mathrm{\Gamma},\rho)$. The representation space is given by 
\begin{equation}
\mathcal{H}_{\mathrm{\Gamma}}=\oplus_{v\in \mathcal{V}}\ \mathcal{H}_v \qquad \text{such that} \qquad \mathcal{H}_v=\ell^2(\mathcal{V})\otimes\mathbb{C}^{2n_v}
\end{equation}
We define the *-representation $\pi_{\mathrm{\Gamma}}$ of $C_{\mathrm{Lip}}(\mathrm{\Gamma})$ on $\mathcal{H}_\mathrm{\Gamma}$ by:
\begin{equation}
\pi (f)\psi(v)=
\oplus_{v_+\sim v}\left( 
\begin{array}{cc}
f(v_+) & 0 \\ 
0 & f(v)
\end{array} 
\right)\psi(v) 
\end{equation}
The Dirac operator $D_\mathrm{\Gamma}$ is an operator on $\mathcal{H}_\mathrm{\Gamma}$ and is given by
\begin{equation}
D_\mathrm{\Gamma}\psi(v)= \oplus_{v_+\sim v}
\frac{1}{\rho(v,v_+)}
\left( 
\begin{array}{cc}
0 & 1 \\ 
1 & 0
\end{array} 
\right)\psi(v)
\end{equation}
The grading operator $\gamma_\mathrm{\Gamma}$ restricted on $\mathcal{H}_v$ is given by 
\begin{equation}
\gamma_\mathrm{\Gamma}|_{\mathcal{H}_v} = 1_{\ell^2(\mathcal{V})}\otimes 
\left( 
\begin{array}{cc}
1 & 0 \\ 
0 & -1
\end{array} 
\right)\otimes 1_{n_v}
\end{equation}
\begin{mprop}
$\pi$ is a faithful $*$-representation of $C_{\mathrm{Lip}}(\mathrm{\Gamma})$.
\end{mprop}
\begin{proof}
It follows immediately that $\pi$ is a $*$-representation for $C_{\mathrm{Lip}}(\mathrm{\Gamma})$. It is bounded, since $f$ is continuous and $\mathrm{\Gamma}$ is finite. Moreover, if $\pi(f)=0$ then $f$ vanishes on the graph $\mathrm{\Gamma}$.
\end{proof}
\begin{mprop}
$(C_{\mathrm{Lip}}(\mathrm{\Gamma}), \mathcal{H}_\mathrm{\Gamma}, D_\mathrm{\Gamma},\gamma_\mathrm{\Gamma})$ is an even spectral triple.
\end{mprop}
\begin{proof}
It follows by definition that $D_{\mathrm{\Gamma}}$ is a self-adjoint operator. We then check that $\gamma_\mathrm{\Gamma}^*=\gamma_\mathrm{\Gamma}$, $\gamma_\mathrm{\Gamma}^2=\gamma_\mathrm{\Gamma}$, $\gamma_{\mathrm{\Gamma}}D_{\mathrm{\Gamma}}=-D_{\mathrm{\Gamma}}\gamma_{\mathrm{\Gamma}}$ and $\mathrm{\Gamma}\pi(f)=\pi(f)\mathrm{\Gamma}$ for all $f\in C_{\mathrm{Lip}}(\mathrm{\Gamma})$.
\end{proof}
\subsection{Inverse limit of spectral triples}
\noindent
Consider two finite $\mathbb{R}$-trees $\mathrm{\Gamma}$ and $\mathrm{\Gamma}'$ such that $\mathrm{\Gamma}\leq \mathrm{\Gamma}'$. As mentioned, there exists a surjective map $r_{\mathrm{\Gamma}\mathrm{\Gamma}'}:\mathrm{\Gamma}'\rightarrow \mathrm{\Gamma}$. Therefore, we can define the pullback map:
\begin{equation}
r^*_{\mathrm{\Gamma}\mathrm{\Gamma}'}:C_{\mathrm{Lip}}(\mathrm{\Gamma})\rightarrow C_{\mathrm{Lip}}(\mathrm{\Gamma}'), \qquad r^*_{\mathrm{\Gamma}\mathrm{\Gamma}'}(f)(x)=f(r_{\mathrm{\Gamma}\mathrm{\Gamma}'}(x))
\end{equation}
which satisfies by construction the following statement.
\begin{mprop} The pullback map $r^*_{\Gamma\Gamma'}$ is injective.
\end{mprop}
\noindent
On the other hand, the inclusion induces a pushforward map $\iota :\mathrm{\Gamma}\rightarrow \mathrm{\Gamma}'$ which is an inclusion of the representation spaces:
\begin{equation}
\iota_{\mathrm{\Gamma}\mathrm{\Gamma}'}:\mathcal{H}_{\mathrm{\Gamma}}\rightarrow \mathcal{H}_{\mathrm{\Gamma}'}, \qquad \iota_{\mathrm{\Gamma}\mathrm{\Gamma}'}(\psi)(x)=
\left\lbrace 
\begin{array}{cc}
\psi(x) & \text{if $x\in \iota_{\mathrm{\Gamma}\mathrm{\Gamma}'}(\mathrm{\Gamma})$,} \\ 
0 & \text{otherwise.}
\end{array} 
\right. 
\end{equation}
We consider the triple $(A_\mathrm{\Gamma},\mathcal{H}_\mathrm{\Gamma},D_\mathrm{\Gamma})$.
\begin{mdef}
A morphism between two spectral triples $(A_i,\mathcal{H}_i,D_i)$ for $i=1,2$ is a pair $(\phi,I)$ consisting of a unital $*$-homomorphism $\phi:A_1\rightarrow A_2$ and a bounded linear operator $I:\mathcal{H}_1\rightarrow \mathcal{H}_2$ satisfying the following conditions:
\begin{itemize}
\setlength\itemsep{0em}
\item[(1)]$\phi(A_1^\infty)\subseteq A_2^\infty$ where $A_1^\infty$ and $A_2^\infty$ are defined in \cite{floricel_inductive_2017};
\item[(2)] $I\pi_1(a)=\pi_2(\phi(a))I$, for every $a\in A_1$;
\item[(3)] $I(Dom(D_1))\subseteq Dom(D_2)$ and $ID_1=D_2I$.
\end{itemize}
A morphism $(\phi,I)$ is said to be isometric if $\phi$ is injective and $I$ is an isometry.\\
A morphism between two even spectral triples $(A_i,\mathcal{H}_i,D_i,\gamma_i)$ for $i=1,2$ is a pair $(\phi,I)$ satisfying the additional condition:
\begin{itemize}
\setlength\itemsep{0em}
\item[(4)] $I\gamma_1=\gamma_2I$;
\end{itemize}
\end{mdef}
\begin{mprop} The pair $(r^*_{\mathrm{\Gamma}\mathrm{\Gamma}'},\iota_{\mathrm{\Gamma}\mathrm{\Gamma}'})$ is a morphism of spectral triples between the two spectral triples $(A_\mathrm{\Gamma},\mathcal{H}_\mathrm{\Gamma},D_\mathrm{\Gamma})$ and $(A_{\mathrm{\Gamma}'},\mathcal{H}_{\mathrm{\Gamma}'},D_{\mathrm{\Gamma}'})$.
\end{mprop}
\begin{proof} We check the different points given in the above definition.
\begin{itemize}
\item[(1)] For any $a\in A_\mathrm{\Gamma}$, the commutator $[D,\pi_\mathrm{\Gamma}(a)]$ is bounded on $\mathcal{H}_\mathrm{\Gamma}$. Similarly, the commutator $[D,\pi_{\mathrm{\Gamma}'}(a)]$ is bounded on $\mathcal{H}_{\mathrm{\Gamma}'}$. Therefore, $A^\infty_\mathrm{\Gamma} = A_\mathrm{\Gamma}$ and $A^\infty_{\mathrm{\Gamma}'} = A_{\mathrm{\Gamma}'}$.
\item[(2)] Using the fact that $\iota_{\mathrm{\Gamma}\mathrm{\Gamma}'}(\psi)(v)\neq 0$ iff $v\in \iota_{\mathrm{\Gamma}\mathrm{\Gamma}'}(\mathrm{\Gamma})$ and the definition of the representation, we check that:
\begin{align*}
&\iota_{\mathrm{\Gamma}\mathrm{\Gamma}'}(\pi_1 (f)\psi)(v)=
\oplus_{\iota_{\mathrm{\Gamma}\mathrm{\Gamma}'}(v_+)\sim \iota_{\mathrm{\Gamma}\mathrm{\Gamma}'}(v)}\left( 
\begin{array}{cc}
f(\iota_{\mathrm{\Gamma}\mathrm{\Gamma}'}(v_+)) & 0 \\ 
0 & f(\iota_{\mathrm{\Gamma}\mathrm{\Gamma}'}(v))
\end{array} 
\right)\psi_{\iota_{\mathrm{\Gamma}\mathrm{\Gamma}'}}(v) \\
&\\
&\pi_2(r^*_{\mathrm{\Gamma}\mathrm{\Gamma}'}(f))\iota_{\mathrm{\Gamma}\mathrm{\Gamma}'}(\psi)(v)=
\oplus_{v_+\sim v}\left( 
\begin{array}{cc}
f\circ r_{\mathrm{\Gamma}\mathrm{\Gamma}'}(v_+) & 0 \\ 
0 & f\circ r_{\mathrm{\Gamma}\mathrm{\Gamma}'}(v)
\end{array} 
\right)\psi_{\iota_{\mathrm{\Gamma}\mathrm{\Gamma}'}}(v)
\end{align*}
and for $v\in \iota_{\mathrm{\Gamma}\mathrm{\Gamma}'}(\mathrm{\Gamma})$, we have $r_{\mathrm{\Gamma}\mathrm{\Gamma}'}(v)=v$.
\item[(3)] Since $\iota_{\mathrm{\Gamma}\mathrm{\Gamma}'}(\ell^2(\mathrm{\Gamma}))\subseteq \ell^2(\mathrm{\Gamma}')$, then $\iota_{\mathrm{\Gamma}\mathrm{\Gamma}'}(Dom(D_1))\subseteq Dom(D_2)$ and
\begin{align*}
\iota_{\mathrm{\Gamma}\mathrm{\Gamma}'}(D\psi)=D \iota_{\mathrm{\Gamma}\mathrm{\Gamma}'}(\psi)
\end{align*}
\noindent
In addition, the map $\iota_{\mathrm{\Gamma}\mathrm{\Gamma}'}$ is an isometry.
\end{itemize}
\end{proof}
\begin{mdef} Let $(J,\leq)$ be a directed index set and $(A_j,\mathcal{H}_j,D_j)_{j\in J}$ be a family of spectral triples. Suppose that for every $j,k\in J$ and $j\leq k$ there exists an isometric morphism $(\phi_{jk},I_{jk})$ from $(A_j,\mathcal{H}_j,D_j)$ to $(A_k,\mathcal{H}_k,D_k)$ satisfying $\phi_{kl}\phi_{jk}=\phi_{jl}$ and $I_{kl}I_{jk}=I_{jl}$, for all $j,k,l\in J$ and $j\leq k\leq l$. The resulting system $\left\lbrace (A_j,\mathcal{H}_j,D_j), (\phi_{jk},I_{jk}) \right\rbrace_{J} $ is called an inductive system of spectral triples.
\end{mdef}
\begin{mdef}
The triple $(A,\mathcal{H},D)$ is called the inductive realization of the inductive system and denoted by $\left\lbrace (A_j,\mathcal{H}_j,D_j), (\phi_{jk},I_{jk}) \right\rbrace_{J}$.
\end{mdef}
\noindent
We consider the inverse limit system of finite trees $(\mathrm{\Gamma}_j,r_{j,j+1})_{j\in J}$ defining the Berkovich projective space $\mathbb{P}^1_{\mathrm{Berk}}$ obtained as the inverse limit 
\begin{equation}
\mathbb{P}^1_{\mathrm{Berk}}=\lim_{\longleftarrow}(\mathrm{\Gamma}_j,r_{j,j+1})_{j\in J}.
\end{equation}
We define the family of spectral triples $\left\lbrace (A_j,\mathcal{H}_j,D_j), (r^*_{jk},\iota_{jk}) \right\rbrace_{J}$ with the following notation:
\begin{equation}
A_j:=C_{\mathrm{Lip}}(\mathrm{\Gamma}_j), \quad \mathcal{H}_j=\ell^2(\mathrm{\Gamma}_j), \quad D_j=D_{\mathrm{\Gamma}_j}
\end{equation}
and with the isometric morphism:
\begin{equation}
r^*_{jk}:A_j\rightarrow A_k, \quad \iota_{jk}: \mathcal{H}\rightarrow \mathcal{H}_k.
\end{equation}
\begin{mprop}
We have the following direct limits:
\begin{equation}
C_\mathrm{Lip}(\mathbb{P}^1_{\mathrm{Berk}}(\mathbb{C}_p)) = \lim_{\longrightarrow}C_\mathrm{Lip}(\mathrm{\Gamma}), \qquad \ell^2(\mathbb{P}^1_\mathrm{Berk})=\lim_{\longrightarrow}\ell^2(\mathrm{\Gamma}), \qquad \pi=\lim_{\longrightarrow}\pi_\mathrm{\Gamma}
\end{equation}
\end{mprop}
\noindent
Therefore, $\pi$ is the unique representation of the $C^*$-algebra $A$ on $\mathcal{H}$ such that
\begin{equation}
\pi(r_j(a))\iota_j =\iota_j\pi_j(a)
\end{equation}
One can describe the inductive limit of the family of operators $\left\lbrace D_j\right\rbrace_{j\in J} $. Following the construction in \cite{floricel_inductive_2017} , we consider the dense domain $\mathscr{D}$ of $\mathcal{H}$,
\begin{equation}
\mathscr{D}=\bigcup_{j\in J}\iota_j(Dom(D_j)).
\end{equation}
For every vector $\psi \in \mathscr{D}$ of the form $\psi=\iota_j\psi_j$, where $\psi_j\in Dom(D_j)$, define 
\begin{equation}
D\psi =\iota_jD_j\psi_j
\end{equation}
It follows that $D$ is a densely defined symmetric operator. Moreover, since the operators $D_j$ are self-adjoint, we have that $Range(D_j\pm i)=\mathcal{H}_j$, for every $j\in J$. Consequently, $Range(D\pm i)$ is dense in $\mathcal{H}$, and thus $D$ is essentially selfadjoint.\par 
\noindent
We shall use the same letter $D$, or the symbol $\displaystyle{D=\lim_{\longrightarrow}D_j}$, to denote the closure of this essentially selfadjoint operator, and call it the inductive limit of the family of operators $\left\lbrace D_j\right\rbrace_{j\in J}$. Similarly, we use the symbol $\gamma$ to denote the inductive limit $ \displaystyle{\gamma=\lim_{\longrightarrow}\gamma_j}$.\par 
\noindent
Putting the various strands together, we consider the following inductive system given by the data $\left\lbrace (A_j,\mathcal{H}_j, D_j,\gamma_j), (\phi_{j,j+1},I_{j,j+1})\right\rbrace_{j\in J}$ and its inductive realization:
\begin{equation*}
(C_{\mathrm{Lip}}(\mathbb{P}^1_\mathrm{Berk}),\ell^2(\mathbb{P}^1_\mathrm{Berk}),\pi, D,\gamma). 
\end{equation*}
\begin{mth}
\label{theorem1}
The operator $D$ is self-adjoint with compact resolvent and $[D,\pi(f)]$ is a bounded operator for any Lipschitz continuous function $a\in C_\mathrm{Lip}(\mathbb{P}^1_\mathrm{Berk})$ i.e. the spectral triple $(C_\mathrm{Lip}(\mathbb{P}^1_\mathrm{Berk}),\ell^2(\mathbb{P}^1_\mathrm{Berk}),\pi,D,\gamma)$ is an even spectral triple.
\end{mth}
\begin{proof}
Let us recall that for an integer $j$, the Dirac operator $D_{\mathrm{\Gamma}_j}$ is defined by
\begin{equation}
D_{\mathrm{\Gamma}_j}\psi(v)= \oplus_{v_+\sim v}
\frac{1}{\rho(v,v_+)}
\left( 
\begin{array}{cc}
0 & 1 \\ 
1 & 0
\end{array} 
\right)\psi(v)
\end{equation}
for any $v \in \mathrm{\Gamma}_j$. Then, we can bound the operator $D_j$ using the metric $\rho$ on the finite graph:
\begin{equation}
\|D_j\|=\sup_{e\in \mathcal{E}_j}\frac{1}{\ell(e)}
\end{equation}
where for an edge $e\in \mathcal{E}_j$ joining the vertices $v_1$ and $v_2$ the length $\ell(e)$ is $\ell(e):=\rho(v_1,v_2)$. We denote by $R_\lambda(D)$ the resolvent of an operator $D$ at a point $\lambda\in \mathbb{C}\backslash\mathbb{R}$. The sequence $\left\lbrace \|D_j\| \right\rbrace_{j\in \mathbb{NN}} $ is unbounded, then using \cite[Thm 3.1]{floricel_inductive_2017} we have that $D$ has compact resolvent if and only if 
\begin{equation*}
\|\iota_j\circ R_\lambda(D_j)\circ \iota_j^*\|=\sup_{k>j}\frac{1}{|\|D_k\|-\lambda|}\rightarrow 0
\end{equation*}
for $\lambda\in \mathbb{C}\backslash\mathbb{R}$, which is equivalent to $\|D_j\|\rightarrow 0$ as $j\rightarrow 0$.\\
Moreover, using the Lipschitz continuity of $a\in C_{\mathrm{Lip}}(\mathrm{\Gamma}_j)$, we have that the family of operators given by $\left\lbrace [D_k,\pi_k(r_{j,k}(a))]\right\rbrace_{k\geq j} $ is uniformly bounded. Then, the operator $D$ is such that $[D,\pi(a)]$ is bounded on a dense subalgebra  of $C_{\mathrm{Lip}}(\mathbb{P}^1_{\mathrm{Berk}}(\mathbb{C}_p))$.
\end{proof}
\section{$\mathbb{P}^1_{\mathrm{Berk}}(\mathbb{C}_p)$ as the universal Wa\.zewski dendrite}
\label{Sect3}
\noindent
 We would like to focus in this work on yet another construction of $\mathbb{P}^1_{\mathrm{Berk}}(\mathbb{C}_p)$ as an inverse limit of compact spaces, called \textit{dendrites}. The present section is based on the work \cite{banic_wazewskis_2013} for which we refer to for proves and more details on the following statements.
 \begin{mdef}
 A \textit{continuum} is a compact, connected metrizable space. A \textit{simple closed curve} in a topological space is any subspace homeomorphic to a circle. A \textit{dendrite} is a locally connected continuum containing no simple closed curve. Dendrites may be thought of as topological generalizations of trees in which branching may occur at a dense set of points.
 \end{mdef}
 \noindent
 The definitions of branch points, ordinary points and endpoints on $\mathbb{R}$-trees, extend to dendrites. We can now introduce the so-called \textit{Wa\.{z}ewski's universal dendrite}.
\begin{mth}[Wa\.{z}ewski's universal dendrite]
Up to homeomorphism, there is a unique dendrite W such that its branch points are dense in W and there are $\aleph_0$ branches at each branch point. The dendrite W is called the Wa\.{z}ewski's universal dendrite.
\end{mth}
\noindent
In particular, the Wa\.{z}ewski's universal dendrite can be obtained from a chain of dendrites $D_1\subseteq D_2 \subseteq D_3 \subseteq \dots$, then defines certain bonding maps $f_n:D_{n+1}\rightarrow D_n$ and then finally obtains Wa\.{z}ewski's universal dendrite as $\lim\limits_{\longleftarrow}\left\lbrace D_n,f_n \right\rbrace_{n=1}^\infty $.
\begin{mprop}
For any dendrite $D$, $D$ is homeomorphic to the Wa\.{z}ewski's universal dendrite if and only if its set of branching points is dense in $D$ and each of its branching points has $\aleph_0$ branches.
\end{mprop}
\noindent
From this, one can immediately notice that $\mathbb{P}^1_{\mathrm{Berk}}(\mathbb{C}_p)$ is a dendrite and the subset of branching points $\mathbb{H}^\mathbb{Q}_{\mathrm{Berk}}(\mathbb{C}_p)$ is dense. Therefore, the Berkovich projective line is homeomorphic to the Wa\.{z}ewski's universal dendrite.\\

\noindent
If $(X,d)$ is a compact metric space, then $2^X$ denotes the set of all nonempty closed subsets of $X$. The set $2^X$ can be equipped with the \textit{Hausdorff metric,} $d_H$ making $(2^X,d_H)$ a metric space called the \textit{hyperspace} of the metric space $(X,d)$.\\
A set valued function from $X$ to $Y$ is a single-valued function $f:X\rightarrow 2^Y$. The \textit{graph} $\mathrm{Gr}(f)$ of a set-valued function $f:X\rightarrow 2^Y$ is the set of all points $(x,y)\in X\times Y$ such that $y\in f(x)$. A function $f:X\rightarrow 2^Y$, where $X$ and $Y$ are compact metric spaces, is an \textit{upper semi-continuous} set-valued function from $X$ to $Y$ (abbreviated \textit{u.s.c.}) if for each open set $V\subseteq Y$, the set $\lbrace x\in X \ | \ f(x)\subseteq V\rbrace$ is an open set in $X$.
\begin{mprop}
Let $X$ and $Y$ be compact metric spaces and $f:X\rightarrow 2^Y$ a set-valued function. Then $f$ is u.s.c. if and only if its graph $\mathrm{Gr}(f)$ is closed in $X\times Y$.
\end{mprop}
\begin{mdef}[Inverse sequence and limit]
An \textit{inverse sequence} of compact metric spaces $X_k$ with u.s.c. bonding functions $f_k$ is a sequence $\left\lbrace X_k,f_k\right\rbrace_{k=1}^\infty$, where $f_k:X_{k+1}\rightarrow 2^{X_k}$.\\
The \textit{inverse limit} of an inverse sequence $\left\lbrace X_k,f_k\right\rbrace_{k=1}^\infty$ with u.s.c. bonding functions is defined as the subspace
\begin{equation*}
\lim_{\longleftarrow}\left\lbrace X_k,f_k\right\rbrace_{k=1}^\infty:=\left\lbrace x=(x_1,x_2,x_3,\dots)\in \prod_{k=1}^\infty X_k : x_k\in f(x_k)\right\rbrace 
\end{equation*}
\end{mdef}
\noindent
In this paper, we look at the special case of inverse limit where $X_k=\left[ 0,1\right] $ and $f_k=f$ for some function $f:[0,1]\rightarrow 2^{[0,1]}$. In such case, the inverse limit  will be denoted by $\lim\limits_{\longleftarrow}\left\lbrace [0,1],f\right\rbrace_{k=1}^\infty$.\\

\noindent
Let $X$ and $Y$ be compact sets, and a given set-valued u.s.c function $f:X\rightarrow 2^Y$ with graph $\mathrm{Gr}(f)$. Let $\pi_1:\mathrm{Gr}(f)\rightarrow X$ and $\pi_2:\mathrm{Gr}(f)\rightarrow Y$ be the coordinate projections. Then, we have 
\begin{equation*}
f(x)=\pi_2(\pi_1^{-1}(x)), \quad x\in X,
\end{equation*} 
Equivalently,  given a closed set $G\subseteq X\times Y$ such that $\pi_1(G)=X$, then $\pi_2\circ \pi_1^{-1}$ defines a set-valued u.s.c. function such that $\mathrm{Gr}(\pi_2\circ \pi_1^{-1})=G$.\\

\noindent
Let $\mathrm{\Delta}$ be the diagonal subset of $[0,1]\times [0,1]$ i.e. the subset
\begin{equation*}
\mathrm{\Delta}=\left\lbrace (t,t)\in [0,1]\times [0,1] : t\in [0,1]\right\rbrace 
\end{equation*}
For a fixed integer $n$, let $\lbrace (a_i,b_i)\rbrace_{i=1}^n$ be a finite sequence in $[0,1]\times [0,1]$, such that $a_i<b_i$ for each $i=1,\dots , n$ and $a_i\neq a_j$ whenever $i\neq j$. Next denote by $L(a_i,b_i)_{i=1}^n$ the union of lines 
\begin{equation*}
L(a_i,b_i)_{i=1}^n=\bigcup_{i=1}^n\left([a_i,b_i]\times \lbrace a_i\rbrace \right)\subseteq [0,1]\times [0,1]
\end{equation*}
Then, we define the set 
\begin{equation}
G(a_i,b_i)_{i=1}^n=\mathrm{\Delta}\cup L(a_i,b_i)_{i=1}^n
\end{equation}
which is closed in the product $[0,1]\times [0,1]$, since it is a union of finitely many closed arcs. Furthermore, $\pi_1(G(a_i,b_i)_{i=1}^n)=\pi_2(G(a_i,b_i)_{i=1}^n)=[0,1]$. Therefore, by the closed graph theorem, there is a surjective u.s.c. function $f_{(a_i,b_i)_{i=1}^n}:[0,1]\rightarrow 2^{[0,1]}$ such that its graph $\mathrm{Gr}(f_{(a_i,b_i)_{i=1}^n})=G(a_i,b_i)_{i=1}^n$.
\begin{mdef}[Comb function]
Let $n$ be a positive integer and $\lbrace (a_i,b_i)\rbrace_{i=1}^n$ be a subset $[0,1]\times [0,1]$, such that $0<a_i<b_i$ for each $i=1,\dots,n$ and $a_i\neq a_j$ whenever $i\neq j$. Then, $f:[0,1]\rightarrow 2^{[0,1]}$ is called an $n$-comb function with respect to $\lbrace (a_i,b_i)\rbrace_{i=1}^n$ if $f=f_{(a_i,b_i)_{i=1}^n}$.
\end{mdef}
\begin{mnot}
For each positive integer $j$, let $i_j$ be a nonnegative integer. We use the notation
\begin{equation*}
(a_1^{i_1},a_2^{i_2},a_3^{i_3},\dots)
\end{equation*}
to represent the point $( \underbrace{a_1,a_1,\dots,a_1}_{i_1}, \underbrace{a_2,a_2,\dots,a_2}_{i_2},\dots) $ and 
\begin{equation*}
(a_1^{i_1},a_2^{i_2},a_3^{i_3},\dots,a_j^{i_j},t^\infty)
\end{equation*}
to denote the point $( \underbrace{a_1,a_1,\dots,a_1}_{i_1}, \underbrace{a_2,a_2,\dots,a_2}_{i_2},\dots,\underbrace{a_j,a_j,\dots,a_j}_{i_j},t,t,t,\dots) $.
\end{mnot}
\begin{mex}
Let $f$ be a 1-comb function with respect to the point $(a,b)$. Then a point $x\in \lim\limits_{\longleftarrow}\lbrace [0,1],f\rbrace_{k=1}^\infty$ if and only if
either $x=(t^\infty)$ for a $t\in [0,1]$, or there exists a positive integer $n$ such that $x=(a^n,t^\infty)$ for a $t\in (a,b]$. Therefore, $\lim\limits_{\longleftarrow}\lbrace [0,1],f\rbrace_{k=1}^\infty$ is a star with center the point $(a^\infty)$ and beams $L_0=\lbrace (t^\infty) \ | \ t\in [0,a]\rbrace$, $L'_0=\lbrace (t^\infty) \ | \ t\in [a,1]\rbrace$ and $L_n=\lbrace (a^n,t^\infty) \ | \ t\in [a,b]\rbrace$. The center $(a^\infty)$ is the only ramification point, and the maximal free arcs are the beams.
\end{mex}
\begin{mprop}
Let $n$ be any positive integer and let $f:[0,1]\rightarrow 2^{[0,1]}$ be any $n$-comb function. Then $\lim\limits_{\longleftarrow}\lbrace [0,1],f\rbrace_{k=1}^\infty$ is a dendrite. In particular, we will denote by $D_n$ the dendrite 
\begin{equation*}
D_n=\lim\limits_{\longleftarrow}\lbrace [0,1],f_{(a_i,b_i)_{i=1}^n}\rbrace_{k=1}^\infty
\end{equation*}
\end{mprop}
\noindent
Let $\lbrace (a_n,b_n)\rbrace_{n=1}^\infty$ be any sequence in $[0,1]\times [0,1]$ such that 
\begin{itemize}
\setlength\itemsep{0em}
\item[1.] $a_n<b_n$ for each $n$
\item[2.] $a_i\neq a_j$ whenever $i\neq j$
\item[3.] $\lim\limits_{n\rightarrow \infty}(b_n-a_n)=0.$
\end{itemize}
Then, similarly one can construct the set $G(a_n,b_n)_{n=1}^\infty$ 
\begin{equation}
G(a_n,b_n)_{n=1}^\infty=\mathrm{\Delta}\cup L(a_n,b_n)_{n=1}^\infty
\label{graph}
\end{equation}
as a closed subset of $[0,1]\times [0,1]$. Therefore, there exists a surjective u.s.c. function $f_{(a_n,b_n)_{n=1}^\infty}:[0,1]\rightarrow 2^{[0,1]}$ such that its graph is exactly $G(a_n,b_n)_{n=1}^\infty$. Thus, $f_{(a_n,b_n)_{n=1}^\infty}$ is called the comb function with respect to the set of points $\lbrace (a_n,b_n)\rbrace_{n=1}^\infty$.
\begin{mth}
Let $f:[0,1]\rightarrow 2^{[0,1]}$ be any comb function with respect to an admissible sequence $\left\lbrace (a_n,b_n)\right\rbrace_{n=1}^\infty $ such that the set $\left\lbrace a_n \ | \ n=1,2,\dots \right\rbrace $ is dense in $[0,1]$. Then,
\begin{equation*}
\lim\limits_{\longleftarrow}\lbrace [0,1],f\rbrace_{k=1}^\infty = \mathrm{Cl}\left( \bigcup_{n=1}^\infty D_n\right)
\end{equation*}
Moreover, $\lim\limits_{\longleftarrow}\lbrace [0,1],f\rbrace_{k=1}^\infty$ is homeomorphic to Wa\.{z}ewski's universal dendrite.
\end{mth}
\begin{mcor}
Let $f:[0,1]\rightarrow 2^{[0,1]}$ be the comb function with respect to a sequence $\lbrace (a_n,b_n)_{n=1}^\infty\rbrace$. Any point $x\in \lim\limits_{\longleftarrow}\lbrace [0,1],f\rbrace_{k=1}^\infty$ can be classified as follows:
\begin{itemize}
\setlength\itemsep{0em}
\item[1.] If $x\in D_n$ for some positive integer $n$ then, either $x=(t^\infty)$, $t\in [0,1]$, or is of the form
\begin{equation*}
x=(a_{i_1}^{k_1},a_{i_2}^{k_2},a_{i_3}^{k_3},\dots,a_{i_m}^{k_m},t^\infty)
\end{equation*}
for some positive integer $m$ and for each $\ell\leq m$ it holds that $i_\ell\leq n$, $k_\ell>0$, $a_{i_\ell}<a_{i_{\ell +1}}\leq b_{i_\ell}$, and $a_{i_m}\leq t \leq b_{i_m}$. In particular,
\begin{itemize}
\setlength\itemsep{0em}
\item[a.] $x$ is a branching point in $D_n$ if and only if 
\begin{equation*}
x=(x_1,x_2,x_3,\dots,x_m,a_j^\infty), \quad \text{$j\leq n$}
\end{equation*}
\item[b.] $x$ is an endpoint in $D_n$ for some positive $n$, i.e.
\begin{equation*}
x=(x_1,x_2,x_3,\dots,x_m,b_i^\infty), \quad \text{$i\leq n$}
\end{equation*}
\end{itemize}
\item[2.] If $x\in \lim\limits_{\longleftarrow}\lbrace [0,1],f\rbrace_{k=1}^\infty\backslash \bigcup_{n=1}^\infty D_n$, then it has the general form
\begin{equation*}
x=(a_{i_1}^{k_1},a_{i_2}^{k_2},a_{i_3}^{k_3},\dots).
\end{equation*}
\end{itemize}
\label{cor1}
\end{mcor}
\begin{mnot}
For any dense subset $\lbrace a_n,n\in \mathbb{N}\rbrace$ of the interval $(0,1)$, there exists a comb function $f$ such that the inverse limit $\lim\limits_{\longleftarrow}\lbrace [0,1],f\rbrace_{k=1}^\infty$ is homeomorphic to the Wa\.zewski's universal dendrite. From now on, we will assume that the comb map $f_{(a_i,b_i)_{i=1}^\infty}$ is constructed such that the set $\lbrace a_n : n\in \mathbb{N}\rbrace$ is an enumeration of the rational numbers in $(0,1)$ i.e. $\mathbb{Q}\cap (0,1)$. In the rest of the paper, we will use $q_n$ instead of $a_n$ to emphasize the fact that we are using rational numbers as coordinates for branch points. 
\end{mnot}
\begin{mrmk}
In Corollary \ref{cor1}, we recognize the general structure of a \textit{dendrite} and the classification of points of the Berkovich line, identified with the Wa\.zewski universal dendrite. We identify the four categories of points:
\begin{itemize}
\setlength\itemsep{0em}
\item[1)]The branching points or Type II points: 
\begin{equation*}
\mathrm{Br}(X)=\left\lbrace x\in X : x=(q_{i_1}^{k_1},q_{i_2}^{k_2},q_{i_3}^{k_3},\dots,q_{i_{m-1}}^{k_{m-1}},q_{i_m}^\infty), \text{\ for each $\ell\leq m,\  q_{i_\ell}<q_{i_{\ell +1}}$}\right\rbrace 
\end{equation*}
\item[2)] The regular points or Type III points:
\begin{equation*}
\mathrm{Reg}(X)=\left\lbrace x\in X : x=(q_{i_1}^{k_1},q_{i_2}^{k_2},q_{i_3}^{k_3},\dots,q_{i_{m}}^{k_{m}},t^\infty), \text{\ $\forall \ell\leq m,\  q_{i_\ell}<q_{i_{\ell +1}}, t\in (q_{i_m},1]\backslash (q_n)$}\right\rbrace
\end{equation*}
\item[3)] The endpoints or Type IV points:
\begin{equation*}
\mathrm{End}(X)=\left\lbrace x\in X : x=(q_{i_1}^{k_1},q_{i_2}^{k_2},q_{i_3}^{k_3},\dots,q_{i_{m}}^{k_{m}},b_j^\infty), \text{\ $\forall \ell\leq m,\  q_{i_\ell}<q_{i_{\ell +1}}, \ q_{i_m}<b_j$}\right\rbrace
\end{equation*}
\item[4)] The limit points or Type I points:
\begin{equation*}
\mathrm{Lim}(X)=\left\lbrace x\in X : x=(q_{i_1}^{k_1},q_{i_2}^{k_2},q_{i_3}^{k_3},\dots,q_{i_{m}}^{k_{m}},\dots), \text{\ $\forall \ell,\  q_{i_\ell}<q_{i_{\ell +1}}$}\right\rbrace 
\end{equation*}
\end{itemize}
\end{mrmk}
\begin{mrmk}
An important property of the Wa\.zewski dendrite $\mathcal{W}$ is the following: it is a homogeneous space in the following sense: the closure of any connected open subset of $\mathcal{W}$ is homeomorphic to $\mathcal{W}$ (see for instance \cite{charatonik_self-homeomorphic_1994}).
\end{mrmk}
\noindent
\begin{mrmk}
The previous results give a description of Berkovich spaces as inverse systems in the compact space $[0,1]$. This construction allows us to introduce the \textit{shift map} induced by the comb function $f$ in the inverse limit. This shift map defines a dynamical system on $\mathbb{P}^1_{\mathrm{Berk}}(\mathbb{C}_p)$.
\subsection{Two equivalent topologies on $\mathbb{P}^1_{\mathrm{Berk}}(\mathbb{C}_p)$}
So far, we have described the topology on the Berkovich line or equivalently the Wa\.zewski dendrite as the \textit{weak topology}. In this short paragraph, we would like to define an equivalent topology given by the so-called cylindrical sets as a subbasis.
\paragraph*{The weak topology on $\mathbb{R}$-trees.} This topology is also known as the \textit{Observer's topology} on a tree $T$. It is the topology generated (in the sense of a subbasis) by sets $U(x,v)$, for $x\in T$ and $v\in T_x$ such that $U(x,v)=\lbrace y\in T : y\neq x, [x,y]\in v\rbrace$. We denote this topology by $\tau_{\mathrm{O}}$.

\paragraph*{Cylindrical topology.} Let $q=q_{i_1}^{k_1},q_{i_2}^{k_2},q_{i_3}^{k_3},\dots,q_{i_{m}}^{k_{m}} $ be a finite word. Consider the \textit{cylinder set} defined by
\begin{equation}
Z(q):=\left\lbrace qx \in \mathbb{P}^1_{\mathrm{Berk}}(\mathbb{C}_p) : x\in \mathbb{P}^1_{\mathrm{Berk}}(\mathbb{C}_p)\right\rbrace 
\end{equation}
Then, define the topology $\tau_{\mathrm{Cyl}}$ generated (as a subbasis) by the cylinder sets.

\begin{mprop}
The topologies $\tau_{\mathrm{O}}$ and $\tau_{\mathrm{Cyl}}$ are equivalent.
\end{mprop}
\begin{proof}
Fix $x,y \in \mathbb{P}^1_\mathrm{Berk}(\mathbb{C}_p)$ and consider the direction set $U(x,[x,y])$. Let $z\in U(x,[x,y])$ such that, without loss of generality, we have $x\leq z \leq y$. By a density argument, we can find $x'$ of Type II such that $x\leq x' \leq z \leq y$. Then, we would have $x'= q_{i_1}^{k_1},q_{i_2}^{k_2},q_{i_3}^{k_3},\dots,q_{i_n}^{\infty} $, $x= q_{i_1}^{k_1},q_{i_2}^{k_2},q_{i_3}^{k_3},\dots q_{i_m}^{k_m},t^{\infty} $, for some $t\leq q_{i_n}$ and also $z=q_{i_1}^{k_1},q_{i_2}^{k_2},q_{i_3}^{k_3},\dots q_{i_n}^{k_n},t^{\infty}$ for some $t\geq q_{i_n}^{k_n}$. Define the word $q=q_{i_1}^{k_1},q_{i_2}^{k_2},q_{i_3}^{k_3},\dots,q_{i_{n}}^{k_n}$ and the cylinder set $Z(q)$. Then, by construction $z\in Z(q)$ and any point $z'\in Z(q)$ is such that $z'\geq x'$ and thus the path $[x,z']$ is equivalent to the path $[x,y]$ as they share the common segment $[x,x']$. Hence, $z\in Z(q)\subseteq U(x,[x,y])$.\\

\noindent
Conversely, let $q=q_{i_1}^{k_1},q_{i_2}^{k_2},q_{i_3}^{k_3},\dots,q_{i_{m}}^{k_{m}} $ and a cylinder set $Z(q)$. Fix $z\in Z(q)$. In addition, define the Type II point given by $x=q_{i_1}^{k_1},q_{i_2}^{k_2},q_{i_3}^{k_3},\dots,q_{i_{m}}^{\infty}$. We next consider $U(x,[x,z])$, the direction of $z$ at $x$. Then, $z\in U(x,[x,z])$. Moreover, note that any $z' \in U(x,[x,z])$ is such that $z'\geq x$ (in the tree's partial order); thus $z'$ can be written $z'=qy$ for some $y\in \mathbb{P}^1_{\mathrm{Berk}}(\mathbb{C}_p)$ and therefore, we have $z'\in Z(q)$. Hence, $ z\in U(x,[x,z]) \subseteq Z(q)$.
\end{proof}
\section{Noncommutative Geometry on $\mathbb{P}^1_{\mathrm{Berk}}(\mathbb{C}_p)$}
\label{Sect4}
There are several ways to associate a $C^*$-algebra to the space $\mathbb{P}^1_{\mathrm{Berk}}(\mathbb{C}_p)$; we will expose some of them that appear to be relevant to the construction of spectral triples. We start with $C^*$-algebras associated to the space of branching points $\mathbb{H}^{\mathbb{Q}}_{\mathrm{Berk}}(\mathbb{C}_p)$ seen as a countable alphabet. We use the construction presented in the works \cite{boava_c-algebras_2023, boava_algebras_2023}.
\end{mrmk}
\subsection{$C^*$-algebras of countable subshifts}
In this paragraph, we will denote by $\mathbb{Q}_1$ the set $\mathbb{Q}\cap (0,1)$. We consider the set of Type 2 points as a countable alphabet. The \textit{shift map} on $\mathbb{Q}_1^{\mathbb{N}}$ is the map $\sigma:\mathbb{Q}_1^\mathbb{N}\rightarrow \mathbb{Q}_1^\mathbb{N}$ given by $\sigma(x_n)=(x_{n+1})$. Now, we notice that $\mathbb{H}^{\mathbb{Q}}_{\mathrm{Berk}}(\mathbb{C}_p)$ is a subset of $\mathbb{Q}_1^{\mathbb{N}}$ and is \textit{invariant} for $\sigma$ i.e. $\sigma(\mathbb{H}^{\mathbb{Q}}_{\mathrm{Berk}}(\mathbb{C}_p))\subseteq \mathbb{H}^{\mathbb{Q}}_{\mathrm{Berk}}(\mathbb{C}_p)$. For an invariant subset $\mathsf{X}\subseteq \mathcal{A}^\mathbb{N}$, we define $\mathcal{L}_n(\mathsf{X})$ as the set of all words of length $n$ that appear in some sequence of $\mathsf{X}$, that is for $\mathbb{H}^{\mathbb{Q}}_{\mathrm{Berk}}(\mathbb{C}_p)$,
\begin{equation*}
\mathcal{L}_n(\mathbb{H}^{\mathbb{Q}}_{\mathrm{Berk}}(\mathbb{C}_p)):=\lbrace (a_0\dots a_{n-1})\in \mathbb{Q}^n:\exists \ x\in \mathbb{H}^{\mathbb{Q}}_{\mathrm{Berk}}(\mathbb{C}_p)\ \text{s.t.}\ (x_0\dots x_{n-1})=(a_0\dots a_{n-1})\rbrace
\end{equation*}
We set $\mathcal{L}_0(\mathbb{H}^{\mathbb{Q}}_{\mathrm{Berk}}(\mathbb{C}_p))=\emptyset$. The \textit{language} of $\mathbb{H}^{\mathbb{Q}}_{\mathrm{Berk}}(\mathbb{C}_p) $ is then the set 
\begin{equation*}
\mathcal{L}_{\mathbb{H}^{\mathbb{Q}}_{\mathrm{Berk}}(\mathbb{C}_p)}:=\bigcup_{n=0}^\infty \mathcal{L}_n(\mathbb{H}^{\mathbb{Q}}_{\mathrm{Berk}}(\mathbb{C}_p))
\end{equation*}
consisting of all finite words that appear in some sequence in $\mathbb{H}^{\mathbb{Q}}_{\mathrm{Berk}}(\mathbb{C}_p)$. Given the subshift $\mathbb{H}^{\mathbb{Q}}_{\mathrm{Berk}}(\mathbb{C}_p)$ over the alphabet $\mathbb{Q}_1$ and $\alpha,\beta \in \mathcal{L}_{\mathbb{H}^{\mathbb{Q}}_{\mathrm{Berk}}(\mathbb{C}_p)}$, we define
\begin{equation*}
C(\alpha,\beta):=\lbrace \beta x\in \mathbb{H}^{\mathbb{Q}}_{\mathrm{Berk}}(\mathbb{C}_p) : \alpha x\in \mathbb{H}^{\mathbb{Q}}_{\mathrm{Berk}}(\mathbb{C}_p)\rbrace
\end{equation*}
In particular, $Z_\beta:=C(\emptyset,\beta)$ is called the \textit{cylinder set} of $\beta$, and $F_\alpha:=C(\alpha,\emptyset)$ the \textit{follower set} of $\alpha$.\\

\noindent
Let $\mathcal{B}_{\mathbb{H}^{\mathbb{Q}}_{\mathrm{Berk}}(\mathbb{C}_p)}$ be the Boolean algebra of subsets of $\mathbb{H}^{\mathbb{Q}}_{\mathrm{Berk}}(\mathbb{C}_p)$ generated by the sets $C(\alpha,\beta)$ for any words $\alpha,\beta\in \mathcal{L}_{\mathbb{H}^{\mathbb{Q}}_{\mathrm{Berk}}(\mathbb{C}_p)}$.
\begin{mdef}
Let $\mathcal{U}$ be a Boolean algebra. The \textit{unital subshift $C^*$-algebra} $\mathcal{O}_\mathsf{X}$ associated with $\mathsf{X}$ is the universal unital $C^*$-algebra generated by projections $\lbrace p_A:A\in \mathcal{U}\rbrace$ and partial isometries $\lbrace s_a:a\in \mathcal{A}\rbrace$ subject to the relations:
\begin{itemize}
\item[(i)] $p_\mathsf{X}=1$, $p_{A\cap B}=p_Ap_B$, $p_{A\cup B}=p_A+p_B-p_{A\cap B}$ and $p_\emptyset =0$, for every $A,B\in \mathcal{U}$;
\item[(ii)] $s_\beta s^*_\alpha s_\alpha s^*_\beta=p_{C(\alpha,\beta)}$ for all $\alpha,\beta \in \mathcal{L}_\mathsf{X}$ for all $\alpha,\beta\in \mathcal{L}_\mathsf{X}$, where $s_\emptyset = 1$ and, for $\alpha=\alpha_1\dots \alpha_n\in \mathcal{L}_\mathsf{X}$, $s_\alpha=s_{\alpha_1}\dots s_{\alpha_n}$.
\end{itemize}
In particular, $s^*_\alpha s_\alpha=p_{C(\alpha,\emptyset)}=p_{F_\alpha}$ and $s_\beta s^*_\beta=p_{C(\emptyset,\beta)}=p_{Z_\beta}$ for all $\alpha,\beta \in \mathcal{L}_\mathsf{X}$.
\end{mdef}
\noindent
We can now define the universal algebra $\mathcal{O}_{\mathbb{H}^{\mathbb{Q}}_{\mathrm{Berk}}(\mathbb{C}_p)}$ associated to the subshift $\mathbb{H}^{\mathbb{Q}}_{\mathrm{Berk}}(\mathbb{C}_p)$, which satisfies the following properties.
\begin{mprop}
\label{prop1}
In $\mathcal{O}_{\mathbb{H}^{\mathbb{Q}}_{\mathrm{Berk}}(\mathbb{C}_p)}$ the following hold:
\begin{itemize}
\item[(i)]$s^*_q s_p=\delta_{q,p}p_{F_q}$ for all $q,p\in \mathbb{Q}_1$
\item[(ii)]$s^*_\alpha s_\alpha$ and $s^*_\beta s_\beta$ commute for all $\alpha, \beta \in \mathcal{L}_{\mathbb{H}^{\mathbb{Q}}_{\mathrm{Berk}}(\mathbb{C}_p)}$
\item[(iii)]$s^*_\alpha s_\alpha$ and $s_\beta s^*_\beta$ commute for all $\alpha, \beta \in \mathcal{L}_{\mathbb{H}^{\mathbb{Q}}_{\mathrm{Berk}}(\mathbb{C}_p)}$
\item[(iv)]$s_\alpha s_\beta=0$ for all $\alpha,\beta \in \mathcal{L}_{\mathbb{H}^{\mathbb{Q}}_{\mathrm{Berk}}(\mathbb{C}_p)}$ such that $\alpha\beta \notin \mathcal{L}_{\mathbb{H}^{\mathbb{Q}}_{\mathrm{Berk}}(\mathbb{C}_p)}$
\item[(v)] $\mathcal{O}_{\mathbb{H}^{\mathbb{Q}}_{\mathrm{Berk}}(\mathbb{C}_p)}$ is generated by the set $\lbrace s_q,s_q^*:q\in \mathbb{Q}_1\rbrace\cup \lbrace 1\rbrace$.
\end{itemize}
\end{mprop}
\begin{proof}
(i) From the definition of a partial isometry, one has $s_q=p_{Z_q}s_q$ and therefore we can write the following equalities $s_q^*s_p=s_q^*p_{Z_q}p_{Z_p}s_p=\delta_{q,p}p_{F_q}$.\\
(ii) Using the point (ii) in the definition, for $\alpha,\beta\in \mathcal{L}_{\mathbb{H}^{\mathbb{Q}}_{\mathrm{Berk}}(\mathbb{C}_p)}$, we have that $s_\alpha^*s_\alpha s^*_\beta s_\beta=p_{F_\alpha}p_{F_\beta}=p_{F_\beta}p_{F_\alpha}$
from which commutativity follows. The proof of (iii) is identical.\\
(iv) Similarly to the proof of (i), we have by definition of a partial isometry that $s_{\alpha}s_{\beta}=s_\alpha p_{F_\alpha}p_{Z_\beta} s_{\beta}=s_\alpha p_{F_\alpha \cap Z_\beta} s_\beta$. The last term in the previous equality vanishes whenever $\alpha\beta \notin \mathcal{L}_{\mathbb{H}^{\mathbb{Q}}_{\mathrm{Berk}}(\mathbb{C}_p)} $.\\
(v) The algebra $\mathcal{O}_{\mathbb{H}^{\mathbb{Q}}_{\mathrm{Berk}}(\mathbb{C}_p)}$ is generated by the partial isometries $s_a$ and the projections $p_A$ for any $A\in \mathcal{B}_{\mathbb{H}^{\mathbb{Q}}_{\mathrm{Berk}}(\mathbb{C}_p)}$. But the projections on the sets $C_{\alpha,\beta}$ can be written as $p_{C(\alpha,\beta)}=s_\beta s_\alpha^*s_\alpha s_\beta^*$. Moreover, the sets $C(\alpha,\beta)$ generate the Boolean algebra, from which the statement follows.
\end{proof}
\begin{mrmk}
The universal $C^*$-algebra $\mathcal{O}_{\mathbb{H}^{\mathbb{Q}}_{\mathrm{Berk}}(\mathbb{C}_p)}$ can be equivalently constructed as an Exel-Laca graph $C^*$-algebra as shown in \cite[Prop. 4.4]{boava_c-algebras_2023}.
\end{mrmk}
\noindent
Moreover, one can associate a Hilbert space to the set $\mathbb{H}^{\mathbb{Q}}_{\mathrm{Berk}}(\mathbb{C}_p)$, which we denote by $\ell^2(\mathbb{H}^{\mathbb{Q}}_{\mathrm{Berk}}(\mathbb{C}_p),dq)$,
where the measure $dq$ is the counting measure on $\mathbb{H}^{\mathbb{Q}}_{\mathrm{Berk}}(\mathbb{C}_p)$. In this section, this Hilbert space will be denoted by $\mathcal{H}$. In particular, the Hilbert space $\mathcal{H}$ admits an orthonormal basis $\lbrace e_x : x\in  \mathbb{H}^{\mathbb{Q}}_{\mathrm{Berk}}(\mathbb{C}_p)\rbrace $ of Kronecker delta functions.
\begin{mprop}
Consider the family of operators $\lbrace S_q\rbrace_{q\in \mathbb{Q}_1}$ and $\lbrace P_Q\rbrace_{Q\in \mathcal{U}}$ on the space $ \ell^2(\mathbb{H}^{\mathbb{Q}}_{\mathrm{Berk}}(\mathbb{C}_p))$ defined by
\begin{equation}
S_q(e_x):=\left\lbrace \begin{array}{cc}
e_{\sigma_q(x)} & \text{if $x\in F_q$} \\ 
0 & \text{otherwise}
\end{array} 
\right. 
\qquad
P_Q(e_x):=\left\lbrace 
\begin{array}{cc}
e_x & \text{if $x\in Q$} \\ 
0 & \text{otherwise}
\end{array} 
\right. 
\end{equation}
where $\lbrace e_x : x \in \mathsf{X}\rbrace$ is an orthonormal basis of $\ell^2(\mathbb{H}^{\mathbb{Q}}_{\mathrm{Berk}}(\mathbb{C}_p))$. Then the map $\pi:\mathcal{O}_{\mathbb{H}^{\mathbb{Q}}_{\mathrm{Berk}}(\mathbb{C}_p)}\rightarrow B(\mathcal{H})$, such that $\pi(s_q)=S_q$ for all $q\in \mathbb{Q}$ and $\pi(p_Q)=P_Q$ for all $Q\in \mathcal{U}$ defines a $^{*}$-representation.
\label{repO}
\end{mprop}
\begin{proof}
We verify immediately that $P_Q$ is a projection satisfying the following properties: $P_{\mathbb{H}^{\mathbb{Q}}_{\mathrm{Berk}}(\mathbb{C}_p)}=1$, $P_{Q\cap R}=P_QP_R$, $P_{Q\cup R}=P_Q + P_R - P_{Q\cap R}$ and $P_\emptyset = 0$, for every $Q,R\in \mathcal{B}_{\mathbb{H}^{\mathbb{Q}}_{\mathrm{Berk}}(\mathbb{C}_p)}$.\\
Moreover, for any $q\in \mathbb{Q}_1$, $S_q$ is a partial isometry between the following initial and final space:
\begin{equation*}
S_q:\overline{\text{span}}\lbrace e_x : x\in F_a\rbrace \rightarrow \overline{\text{span}}\lbrace e_x : x\in Z_a\rbrace 
\end{equation*}
with the adjoint given by
\begin{equation*}
S_q^*(e_x)=\left\lbrace  
\begin{array}{cc}
e_y, & \text{if $x=ay\in Z_a$} \\ 
0, & \text{otherwise}.
\end{array} 
\right. 
\end{equation*}
In addition, we can define $S_\alpha:=S_{\alpha_1}\cdots S_{\alpha_n}$ for $\alpha=\alpha_1\cdots \alpha_n\in \mathcal{L}_{\mathbb{H}^{\mathbb{Q}}_{\mathrm{Berk}}(\mathbb{C}_p)}$. The partial isometries also generate the projections on the generating sets $S_\beta S_\alpha^*S_\alpha S_\beta^*=P_{C(\alpha,\beta)}$. The claim then follows from the universal property of $\mathcal{O}_{\mathbb{H}^{\mathbb{Q}}_{\mathrm{Berk}}(\mathbb{C}_p)}$.
\end{proof}
\subsection{Semibranching systems on $\mathbb{P}^1_{\mathrm{Berk}}(\mathbb{C}_p)$}
\label{Sect5}
We can extend the previous representation on the set of rational points $\mathbb{H}_{\mathrm{Berk}}^\mathbb{Q}(\mathbb{C}_p)$ to the full space $\mathbb{P}^1_{\mathrm{Berk}}(\mathbb{C}_p)$. We can also exploit the metric structure on the Berkovich line in order to single out invariant measures. To exhibit representation spaces of such $C^*$-algebras, we rely on semibranching systems; we refer to \cite{bratteli_iterated_1996} for further details on such construction.
\subsubsection{Invariant measure}
In order to define a representation of the universal $C^*$-algebra $\mathcal{O}_{\mathbb{H}_{\mathrm{Berk}}^\mathbb{Q}(\mathbb{C}_p)}$ using the space $\mathbb{P}^1_{\mathrm{Berk}}(\mathbb{C}_p)$, we need to single out a measure $\mu$ on the projective line. In order to do so, we will use the isometry group $\mathrm{Iso}(\mathbb{P}^1_{\mathrm{Berk}}(\mathbb{C}_p))$.\\

\noindent
A fact of fundamental importance is that the action of a non-constant rational map on $\mathbb{P}^1(\mathbb{C}_p)$ extends naturally to an action on $\mathbb{P}^1_{\mathrm{Berk}}(\mathbb{C}_p)$, and such map will preserve the type of the point upon which it acts. \\
 
\noindent
Let $\varphi\in \mathbb{C}_p(T)$ be a rational function of degree $d\geq 1$. The usual action on $\mathbb{P}^1(\mathbb{C}_p)$ extends to a continuous action $\varphi:\mathbb{P}^1_{\mathrm{Berk}}(\mathbb{C}_p)\rightarrow \mathbb{P}^1_{\mathrm{Berk}}(\mathbb{C}_p)$.\\

\noindent
In particular, the maps $\gamma \in \mathrm{PGL}_2(\mathbb{C}_p)$ act transitively on type II points of $\mathbb{P}^1_{\mathrm{Berk}}(\mathbb{C}_p)$, and any type II point $\zeta_{a,r}$ can be written as $\zeta_{a,r}=\gamma(\zeta_{root})$, where $\gamma=\left( \begin{array}{cc}
q & a \\ 
0 & 1
\end{array}\right)  $ and $|q|_v=r$.
\begin{mcor}[\cite{baker_potential_nodate}]
$\mathrm{Aut}(\mathbb{P}^1_{\mathrm{Berk}}(\mathbb{C}_p))\simeq \mathrm{PGL}_2(\mathbb{C}_p)$, the group of Möbius transformations (or linear fractional transformations) acting on $\mathbb{P}^1_{\mathrm{Berk}}(\mathbb{C}_p)$.
\end{mcor}
\noindent
The group $\mathrm{PGL}_2(\mathbb{C}_p)$ of Möbius transformations acts continuously on $\mathbb{P}^1_{\mathrm{Berk}}(\mathbb{C}_p)$ in a natural way compatible with the usual action $\mathbb{P}^1(\mathbb{C}_p)$, and this actions preserves $\mathbb{H}_{\mathrm{Berk}}(\mathbb{C}_p)$. Using the definition of $\mathbb{P}^1_{\mathrm{Berk}}(\mathbb{C}_p)$ in terms of multiplicative seminorms (and extending each $[\cdot]_x$ to a seminorm on its local ring in the quotient field $\mathbb{C}_p(T)$) we have 
\begin{equation}
[f]_{M(x)}=[f\circ M]_x
\end{equation}
for each $M\in \mathrm{PGL}_2(\mathbb{C}_p)$. The action of $\mathrm{PGL}_2(\mathbb{C}_p)$ on $\mathbb{P}^1_{\mathrm{Berk}}(\mathbb{C}_p)$ can also be described concretely in terms of Berkovich’s classification theorem, using the fact that each $M\in \mathrm{PGL}_2(\mathbb{C}_p)$ takes closed discs to closed discs. An important observation is that $\mathrm{PGL}_2(\mathbb{C}_p)$ acts \textit{isometrically} on $\mathbb{H}_{\mathrm{Berk}}(\mathbb{C}_p)$ i.e.
\begin{equation}
\rho(Mx,My)=\rho(x,y)
\end{equation}
for all $x,y\in \mathbb{H}_{\mathrm{Berk}}(\mathbb{C}_p)$ and all $M\in \mathrm{PGL}_2(\mathbb{C}_p)$, where $\rho$ is the \textit{big metric}. This shows that the path metric $\rho$ is "coordinate-free".\\

\noindent
Following \cite{fremlin_measure_2006}, since $(\mathbb{P}^1_{\mathrm{Berk}}(\mathbb{C}_p),\rho)$ is a compact metric space with isometry group $\mathrm{PGL}_2(\mathbb{C}_p)$, there exists a non-zero $\mathrm{PGL}_2(\mathbb{C}_p)$-invariant Radon measure on $\mathbb{P}^1_{\mathrm{Berk}}(\mathbb{C}_p)$, denoted by $\mu_G$. Moreover, since $\mathrm{PGL}_2(\mathbb{C}_p)$ acts transitively, then $\mu_G$ is strictly positive. We can then define the Hilbert space $L^2(\mathbb{P}^1_{\mathrm{Berk}}(\mathbb{C}_p),\mu_G)$.
\subsubsection{Representation in $L^2(\mathbb{P}^1_{\mathrm{Berk}}(\mathbb{C}_p),\mu_G)$}
Consider the shift map $\sigma:\mathbb{P}^1_{\mathrm{Berk}}(\mathbb{C}_p)\rightarrow \mathbb{P}^1_{\mathrm{Berk}}(\mathbb{C}_p)$; it is a continuous surjection. Keeping the notations introduced in the previous section, we consider the sets
\begin{equation*}
C(\alpha,\beta):=\lbrace \beta x\in \mathbb{P}^1_{\mathrm{Berk}}(\mathbb{C}_p) : \alpha x\in \mathbb{P}^1_{\mathrm{Berk}}(\mathbb{C}_p)\rbrace
\end{equation*}
for every $\alpha,\beta \in \mathcal{L}_{\mathbb{H}^{\mathbb{Q}}_{\mathrm{Berk}}(\mathbb{C}_p)}$. In particular, $Z_\beta:=C(\emptyset,\beta)$ is the \textit{cylinder set} of $\beta$, and $F_\alpha:=C(\alpha,\emptyset)$ the \textit{follower set} of $\alpha$. Let us define the maps $\sigma_q$ for $q\in $such that
\begin{equation}
\sigma_q:F(q)\rightarrow Z(q), \quad \sigma_q(x)=qx.
\end{equation}
Let $\mathcal{B}_{\mathbb{P}^1_{\mathrm{Berk}}(\mathbb{C}_p)}$ be the Boolean algebra of subsets of $\mathbb{P}^1_{\mathrm{Berk}}(\mathbb{C}_p)$ generated by all the sets $C(\alpha,\beta)$ for $\alpha,\beta\in \mathcal{L}_{\mathbb{H}^{\mathbb{Q}}_{\mathrm{Berk}}(\mathbb{C}_p)}$. similarly, one can construct the universal $C^*$-algebra $\mathcal{O}_{\mathbb{P}^1_{\mathrm{Berk}}(\mathbb{C}_p)}$.
\begin{mth}
\label{mth2}
The $C^*$-algebras $\mathcal{O}_{\mathbb{P}^1_{\mathrm{Berk}}(\mathbb{C}_p)}$ and $\mathcal{O}_{\mathbb{H}^{\mathbb{Q}}_{\mathrm{Berk}}(\mathbb{C}_p)}$ are *-isomorphic.
\end{mth} 
\begin{proof}
The map $j$ defined on cylinder sets such that 
\begin{equation*}
j(\lbrace \beta x\in \mathbb{H}^{\mathbb{Q}}_{\mathrm{Berk}}(\mathbb{C}_p) : \alpha x\in \mathbb{H}^{\mathbb{Q}}_{\mathrm{Berk}}(\mathbb{C}_p)\rbrace)=\lbrace \beta x\in \mathbb{P}^1_{\mathrm{Berk}}(\mathbb{C}_p) : \alpha x\in \mathbb{P}^1_{\mathrm{Berk}}(\mathbb{C}_p)\rbrace
\end{equation*}
extends to an isomorphism of Boolean algebras between $\mathcal{B}_{\mathbb{H}^{\mathbb{Q}}_{\mathrm{Berk}}(\mathbb{C}_p)}$ and $\mathcal{B}_{\mathbb{P}^1_{\mathrm{Berk}}(\mathbb{C}_p)}$. The statement follows then from the universal property of $\mathcal{O}_{\mathbb{H}^{\mathbb{Q}}_{\mathrm{Berk}}(\mathbb{C}_p)}$.
\end{proof}
\noindent
From now on, we will consider $\mathcal{O}_{\mathbb{P}^1_{\mathrm{Berk}}(\mathbb{C}_p)}$ as the $C^*$-algebra associated to the Berkovich line. We will now study its representation on $L^2(\mathbb{P}^1_{\mathrm{Berk}}(\mathbb{C}_p),\mu_G)$.
\begin{mdef}[Semibranch system]
Consider a measure space $(X,\mu)$ and a countable family $\lbrace\sigma_i\rbrace_{i\in \mathbb{N}}$, of measurable maps $\sigma_i:D_i\rightarrow X$, defined on measurable subsets $D_i\subset X$. The family $\lbrace\sigma_i\rbrace_{i\in \mathbb{N}}$ is called \textit{a semibranching system} if the following holds 
\begin{itemize}
\item[(1)] There exists a corresponding family $\lbrace R_i\rbrace_{i\in \mathbb{N}}$  of measurable subsets of $X$ with the property,
\begin{equation}
\mu(X\backslash \cup_i R_i)=0, \quad \text{and} \quad \mu(R_i\cap R_j)=0, \quad \text{for $i\neq j$}
\end{equation}
where we denote by $R_i$ the range $R_i=\sigma_i(D_i)$.
\item[(2)] There is a Radon-Nikodym derivative
\begin{equation*}
\mathrm{\Phi}_{\sigma_i}=\frac{d(\mu\circ \sigma_i)}{d\mu}
\end{equation*}
with $\mathrm{\Phi}_{\sigma_i}>0$, $\mu$-almost everywhere on $D_i$.
\end{itemize}
A measurable map $\sigma:X\rightarrow X$ is called a coding map for the family $\lbrace \sigma_i\rbrace$ if for all $x\in D_i$, $\sigma\circ \sigma_i(x)=x$ for all $x\in D_i$.
\end{mdef}
\begin{mlem}
\label{lemmaRN}
For every $q\in \mathbb{Q}_1$, the Radon-Nikodym derivative $\Phi_{\sigma_q}=1$.
\end{mlem}
\begin{proof}
The sets $F(q)$ and $D(q)$ are both non-empty compact, locally connected metric spaces with no closed curves, hence both are dendrites. Moreover, their sets of branching points are dense (by construction) with infinite branching at each point. By universality of the Wa\.zewski dendrite $\mathcal{W}$, $F(q)$ and $D(q)$ are both homeomorphic to $\mathcal{W}$. In, addition the map $\sigma_q:F(q)\rightarrow D(q)$ is continuous, surjective and thus a homeomorphism, with inverse $\sigma$ (see \cite{charatonik_self-homeomorphic_1994}). Thus, $\sigma_q\in \mathrm{Homeo}(\mathcal{W})$ and by identification with $\mathbb{P}^1_{\mathrm{Berk}}(\mathbb{C}_p)$, $\sigma_g$ can be identified with an element $\gamma_g$ of $\mathrm{PGL}_2(\mathbb{C}_p)$. Using the invariance of $\mu_G$ we can compute the Radon-Nikodym derivative:
\begin{equation}
\Phi_{\sigma_q}=\frac{d(\mu_G\circ \sigma_q)}{d\mu_G}=\frac{d(\mu_G\circ \gamma_g)}{d\mu_G}=1.
\end{equation}
\end{proof}
\begin{mprop}
The shift map $\sigma:\mathbb{P}^1_{\mathrm{Berk}}(\mathbb{C}_p)\rightarrow \mathbb{P}^1_{\text{Berk}}(\mathbb{C}_p)$ is the coding map of the semibranching function system given by the family of maps $\lbrace \sigma_{q}\rbrace_{q}$ such that
\begin{equation}
\sigma_{q}:F(q)\rightarrow Z(q), \qquad \sigma_q(q_1,q_2,\dots)=(q,q_1,q_2,\dots) \qquad \text{for all $q \in \mathbb{Q}_1$}
\end{equation}
Moreover, one can construct a family of operator $(S_{q})_{q\in \mathbb{Q}}$ acting on the Hilbert space $L^2(\mathbb{P}^1_{\mathrm{Berk}}(\mathbb{C}_p),\mu_G)$:
\begin{equation}
(S_{q}\psi)(x)=\chi_{Z(q)}(x)\psi\circ\sigma(x)
\end{equation}
We will also denote the space $L^2(\mathbb{P}^1_{\mathrm{Berk}}(\mathbb{C}_p),\mu_G)$ by $\mathcal{H}$.
\end{mprop}
\begin{proof}
Item (1) in the definition of a semibranching system follows from the fact that the collection $\lbrace Z(q) \ : \ q\in \mathbb{Q}_1 \rbrace$ forms a disjoint partition of $\mathbb{P}^1_{\mathrm{Berk}}(\mathbb{C}_p)$. Item (2) follows from Lemma \ref{lemmaRN}.
\end{proof}
\begin{mlem}
The adjoint of the operator $S_q$ is given by
\begin{equation}
S_{q}^*:\mathcal{H}\rightarrow \mathcal{H}, \quad (S_{q}^*\varphi)(x)=\chi_{F(q)}(x)\psi\circ\sigma_{q}(x)
\end{equation}
\end{mlem}
\begin{proof}
For any $q\in \mathbb{Q}_1$, we can write 
\begin{align*}
\left\langle S_q\psi,\varphi  \right\rangle = \int_{Z(q)}\psi(\sigma(x))\varphi(x)d\mu_G(x)=  \int_{F(q)}\psi(u)\varphi(\sigma_q(x))\mathrm{\Phi}_{\sigma_q}(u) d\mu_G(u) = \left\langle \psi,S_q^*\varphi  \right\rangle
\end{align*}
where we used Lemma \ref{lemmaRN} to state that $\mathrm{\Phi}_{\sigma_q}(u) =1$.
\end{proof}
\noindent
Consider then the $C^*$-algebra generated by the set of partial isometries $\lbrace S_{q}, S_{q}^* : \text{for all $\zeta\leq \zeta'$}\rbrace$ and then define the concrete $C^*$-algebra 
\begin{equation}
A(\mathcal{H}):=\overline{\lbrace S_{q}, S_{q}^* : \text{for all $q\in \mathbb{Q}_1$}\rbrace}^{\mathrm{S.O.T}}
\end{equation}
\begin{mprop}
The operators $S_{q}$ and their adjoints satisfy the relation $S_{q}S_{q}^*=P_{Z(q)}$, where $P_{Z(q)}$ is the projection given by multiplication by the characteristic function $\chi_{Z(q)}$. This gives
\begin{equation}
\sum_{q}S_{q}S_{q}^*=1
\end{equation}
Similarly, $S^*_{q}S_{q}=P_{F(q)}$, where $P_{F(q)}$ is the projection given by multiplication by $\chi_{F(q)}$.
\end{mprop}
\begin{proof}
For any $q\in \mathbb{Q}_1$, we can write 
\begin{equation*}
S_qS_q^*\psi(x)=\chi_{Z(q)}(x)\chi_{F(q)}(\sigma(x))\psi(\sigma(\sigma_q(x)))=\chi_{Z(q)}(x)\chi_{F(q)}(\sigma(x))\psi(x)
\end{equation*}
where we have used the fact that $\sigma (\sigma_q(x))=x$ by definition of $\sigma_q$. Now, if $x\in Z(q)$ then $\sigma(x)\in F(q)$ and thus $\chi_{F(q)}(\sigma(x))=1$. Hence, we deduce that
\begin{equation*}
S_qS_q^*\psi(x)=\chi_{Z(q)}(x)\psi(x)
\end{equation*}
and therefore $S_qS_q^*=P_{Z(q)}$. The proof that $S_q^*S_q=P_{F(q)}$ is identical. Finally, for any $\psi \in \mathcal{H}$ and $x\in \mathbb{P}^1_{\mathrm{Berk}}(\mathbb{C}_p)$, we have 
\begin{equation*}
\sum_{q\in \mathbb{Q}_1}S_qS_q^*\psi(x)=\sum_{q\in \mathbb{Q}_1}P_{Z(q)}\psi(x)=\psi(x)
\end{equation*}
since $\lbrace Z(q) \ : \ q\in \mathbb{Q}_1 \rbrace$ forms a partition of $\mathbb{P}^1_{\mathrm{Berk}}(\mathbb{C}_p)$.
\end{proof}
\begin{mcor}
The map $\pi:\mathcal{O}_{\mathbb{P}^1_{\mathrm{Berk}}(\mathbb{C}_p)}\rightarrow A(\mathcal{H})$ given by 
\begin{equation}
\pi(s_q)=S_q
\end{equation}
 is a *-homomorphism and defines a representation of $\mathcal{O}_{\mathbb{P}^1_{\mathrm{Berk}}(\mathbb{C}_p)}$ as a concrete $C^*$-subalgebra of the algebra of bounded operators $B(L^2(\mathbb{P}^1_{\mathrm{Berk}}(\mathbb{C}_p),\mu_G))$.
\end{mcor}
\subsection{Projection-valued measures}
\noindent
In this section, we denote by $\mathcal{H}$ the Hilbert space $L^2(\mathbb{P}^1_{\mathrm{Berk}}(\mathbb{C}_p),\mu_G)$.
\begin{mdef}
Let $\mathcal{H}$ be a Hilbert space and 
\begin{equation*}
\mathcal{P}_\mathcal{H}:=\left\lbrace P\in B(\mathcal{H}):P=P^2=P^*\right\rbrace 
\end{equation*}
be the set of all orthogonal projections on $\mathcal{H}$. Further, let $\mathrm{\Sigma}(X)$ be the $\sigma$-algebra of a measurable space $X$. An operator-valued map $P:\mathrm{\Sigma}(X)\rightarrow \mathcal{P}_\mathcal{H}$ defined on $\mathrm{\Sigma}(X)$ with values in bounded linear operators on Hilbert space $\mathcal{H}$ is called a \textit{spectral measure} or a \textit{projection valued measure} if 
\begin{itemize}
\item[(1)] $P(X)=1$ and $P(\emptyset)=0$,
\item[(2)] If $B_1,B_2,\dots$ in $\mathrm{\Sigma}(X)$, such that $B_i\cap B_j=\emptyset$ for $i\neq j$, one has
\begin{equation}
P\left( \bigcup_{i=1}^\infty B_i\right)=\sum_{i=1}^\infty P(B_i) 
 \end{equation}
in the strong topology sense.
\item[(3)] $P(E\cap F)=P(E)P(F)$ for $E,F\in \mathrm{\Sigma}(X)$.
\end{itemize}
\end{mdef}
\noindent
As an example, consider $(\mathbb{P}^1_{\mathrm{Berk}}(\mathbb{C}_p),\mu_G)$ as a $\sigma$-finite measure space and let $\mathrm{\Sigma}_G$ be the $\sigma$-algebra of $\mu_G$-measurable subsets of $\mathbb{P}^1_{\mathrm{Berk}}(\mathbb{C}_p)$. For each measurable set $E\subseteq \mathbb{P}^1_{\mathrm{Berk}}$, define $P_G(E)$ the projection in $B(\mathcal{H})$ given by $P(E)=m_{\chi_E}$, that is 
\begin{equation}
P_G(E)f=\chi_Ef.
\end{equation}
Then $E\rightarrow P_G(E)$ is a projection-valued measure on $\mathcal{H}$ and called the \textit{canonical projection-valued measure} on $L^2(\mu_G)$.\\

\noindent
Let $\mathrm{\Sigma}(\mathbb{P}^1_{\mathrm{Berk}}(\mathbb{C}_p))$ be the \textit{cylinder $\sigma$-algebra,} i.e. the $\sigma$-algebra generated by the cylinder sets.
\begin{mprop}
The operator-valued map defined on the cylinders $Z(q)\mapsto S_qS_q^*$ extends to a spectral measure $P:\mathrm{\Sigma}(\mathbb{P}^1_{\mathrm{Berk}}(\mathbb{C}_p))\rightarrow \mathcal{P}_{\mathcal{H}}$. Moreover, we have the equality
\begin{equation}
P(E)=P_G(E ), \qquad \forall E\in \mathrm{\Sigma}(\mathbb{P}^1_{\mathrm{Berk}}(\mathbb{C}_p)).
\end{equation}
\end{mprop}
\begin{proof}
Let $f\in \mathcal{H}$, we use Kolmogorov extension theorem on the map $E\mapsto \left\langle P(E)f,f \right\rangle $ defined on cylinder sets. We just have to check the consistency conditions. Consider two different words $Q=q_1\dots q_n$ and $Q'=q_1'\dots q_n'$ in $\mathcal{L}_{\mathbb{H}^{\mathbb{Q}}_{\mathrm{Berk}}(\mathbb{C}_p)}$, then by definition the cylinders $Z(Q)$ and $Z(Q')$ are disjoint. The projections $P_{Z(Q)}=S_QS_Q^*$ and $P_{Z(Q')}=S_{Q'}S_{Q'}^*$ are orthogonal.\\

\noindent
Again, by construction, for a finite word $Q=q_1\dots q_n$, the cylinder $Z(Q)$ is the disjoint union of the cylinders $Z(q_1\dots q_n q)$ over $q\geq q_n$. We can then write, 
\begin{equation}
\left\langle \sum_{q\geq q_n}P(Z(q_1\dots q_n q))f,f\right\rangle = \left\langle\sum_{q\geq q_n}  S_{Qq}S_{Qq}^*f,f\right\rangle =\left\langle S_QS_Q^*f,f\right\rangle =\left\langle P(Z(Q))f,f\right\rangle 
\end{equation}
Therefore, the consistency relations are true for an arbitrary $f\in \mathcal{H}$ and hence $P$ extends to a projection  valued measure on $\mathcal{P}_{\mathcal{H}}$.\\

\noindent
Finally for $f\in \mathcal{H}$ and for any $E\in \mathrm{\Sigma}(\mathbb{P}^1_{\mathrm{Berk}}(\mathbb{C}_p)) $, using the representation relations, we have $P(E)f=\chi_Ef$ and thus $P$ coincides with $P_G$.
\end{proof}
\begin{mlem}
There exists a \textit{cyclic vector} $f\in \mathcal{H}$ such that 
\begin{equation}
\mu_G(E)=P^f_G(E):=\left\langle P_G(E)f,f \right\rangle =\|P_G(E)f\|^2 
\end{equation}
for every $\mu_G$-measurable sets in $\mathbb{P}^1_{\mathrm{Berk}}$.
\end{mlem}
\begin{proof}
Let us recall that a vector $f\in \mathcal{H}$ is called a \textit{cyclic} vector for the operator $P_G$ if the linear span of the vectors $P_G(E)f$, for $E\in \mathrm{\Sigma}(\mathbb{P}^1_{\mathrm{Berk}}(\mathbb{C}_p)) $, is dense in $\mathcal{H}$. Since $\mu_G$ is a finite measure and the set of simple functions is dense in $\mathcal{H}$, then the identity function $1$ is a cyclic vector for $P_G$.
\end{proof}
\begin{mprop}
Consider the projection-valued measure $(P,\mathrm{\Sigma}(\mathbb{P}^1_{\mathrm{Berk}}(\mathbb{C}_p)))$ with a cyclic vector $f\in \mathcal{H}$. The measure $\mu_G$ coincides withx the real valued Borel measure on the space $\mathbb{P}^1_{\mathrm{Berk}}(\mathbb{C}_p)$ defined by 
\begin{equation}
\mu_f(E):=P^f(E):=\left\langle P(E)f,f \right\rangle =\|P(E)f\|^2 
\end{equation}
Moreover, the measure $\mu_f$ satisfies
\begin{equation*}
\int_{\mathbb{P}^1_{\mathrm{Berk}}(\mathbb{C}_p)}\psi \ d\mu_f= \int_{\mathbb{P}^1_{\mathrm{Berk}}(\mathbb{C}_p)}\psi \ d\mu_G \quad \text{and,} \quad \sum_{q\in \mathbb{Q}_1}\int_{\mathbb{P}^1_{\mathrm{Berk}}(\mathbb{C}_p)}\psi\circ \sigma_q \ d\mu_{S_q^*f}=\int_{\mathbb{P}^1_{\mathrm{Berk}}(\mathbb{C}_p)}\psi \ d\mu_f
\end{equation*}
\end{mprop}
\begin{proof}
Let us denote by $\pi:C(\mathbb{P}^1_{\mathrm{Berk}}(\mathbb{C}_p))\rightarrow B(\mathcal{H})$ the left-multiplication representation of continuous functions on $\mathbb{P}^1_{\mathrm{Berk}}(\mathbb{C}_p)$. We can write 
\begin{align*}
\sum_{q\in \mathbb{Q}_1}\int_{\mathbb{P}^1_{\mathrm{Berk}}(\mathbb{C}_p)}\psi\circ \sigma_q \ d\mu_{S_q^*f} &= \sum_{q\in \mathbb{Q}_1}\int_{\mathbb{P}^1_{\mathrm{Berk}}(\mathbb{C}_p)}\chi_{F(q)}\psi\circ \sigma_q \ d\mu_{S_q^*f}\\ 
&= \sum_{q\in \mathbb{Q}_1}\left\langle S_q^*f, P(F(q))\pi(\psi\circ \sigma_q)S_q^*f\right\rangle 
\end{align*}
Moreover, using the fact that $P(F(q))=S_q^*S_q$, we can say that 
\begin{equation*}
P(F(q))\pi(\psi\circ \sigma_q)S_q^*f=\chi_{F(q)}\pi(\psi\circ \sigma_q)(f\circ \sigma_q)=S_q^*\pi(\psi)f
\end{equation*}
and thus $\left\langle S_q^*f, P(F(q))\pi(\psi\circ \sigma_q)S_q^*f\right\rangle = \left\langle f,S_qS_q^*\pi(\psi)f \right\rangle = \left\langle f,P(Z(q))\pi(\psi)f \right\rangle $. Finally, we conclude that 
\begin{equation*}
\sum_{q\in \mathbb{Q}_1}\int_{\mathbb{P}^1_{\mathrm{Berk}}(\mathbb{C}_p)}\psi\circ \sigma_q \ d\mu_{S_q^*f} = \sum_{q\in \mathbb{Q}_1} \left\langle f,P(Z(q))\pi(\psi)f \right\rangle =\left\langle f,\pi(\psi)f\right\rangle 
\end{equation*}
In other words, we can write 
\begin{equation*}
\sum_{q\in \mathbb{Q}_1}\int_{\mathbb{P}^1_{\mathrm{Berk}}(\mathbb{C}_p)}\psi\circ \sigma_q \ d\mu_{S_q^*f}=\int_{\mathbb{P}^1_{\mathrm{Berk}}(\mathbb{C}_p)}\psi \ d\mu_f
\end{equation*}
\end{proof}
\begin{mprop}
There exists a unique continuous linear map
\begin{equation*}
\widehat{P}:L^\infty(\mathbb{P}^1_{\mathrm{Berk}}(\mathbb{C}_p))\rightarrow B(\mathcal{H})
\end{equation*}
with $\widehat{P}(\chi_E)=P(E)$ for $E\in \mathrm{\Sigma}(\mathbb{P}^1_{\mathrm{Berk}}(\mathbb{C}_p))$. This map is called the spectral integral, and we also write it:
\begin{equation}
\widehat{P}(f)=\int_{\mathbb{P}^1_{\mathrm{Berk}}(\mathbb{C}_p)}f(\zeta)dP(\zeta)
\end{equation}
This map satisfies $\widehat{P}(f)^*=\widehat{P}(\overline{f})$, $\widehat{P}(fg)=\widehat{P}(f)\widehat{P}(g)$ and $\|\widehat{P}(f)\|\leq \|f\|_\infty$ for $f,g\in L^\infty(\mathbb{P}^1_{\mathrm{Berk}}(\mathbb{C}_p))$. In particular, $(\widehat{P},\mathcal{H})$ is a representation of the $C^*$-algebra $L^\infty(\mathbb{P}^1_{\mathrm{Berk}}(\mathbb{C}_p))$ of all bounded measurable functions on the space $\mathbb{P}^1_{\mathrm{Berk}}(\mathbb{C}_p)$.
\end{mprop}
\begin{proof}
The properties of the map $\widehat{P}$ are obtained from the properties of multiplication operators by functions in $L^\infty(\mathbb{P}^1_{\mathrm{Berk}}(\mathbb{C}_p))$. Uniqueness of the spectral integral follows from the density of the subspace of step functions $\left\lbrace \chi_E \ : \ E\in \mathrm{\Sigma}(\mathbb{P}^1_{\mathrm{Berk}}(\mathbb{C}_p)) \right\rbrace $ in $L^\infty(\mathbb{P}^1_{\mathrm{Berk}}(\mathbb{C}_p))$.
\end{proof}
\begin{mcor}
Consider the projection-valued measure $(P,\mathrm{\Sigma}(\mathbb{P}^1_{\mathrm{Berk}}(\mathbb{C}_p)))$ with a cyclic vector $f\in \mathcal{H}$. Define the Hilbert space $\mathcal{H}_f$ such that
\begin{equation}
\mathcal{H}_f:=\overline{\mathrm{span}}\left\lbrace S_\alpha S_\alpha^*f : \alpha\in \mathcal{\mathbb{H}^{\mathbb{Q}}_{\mathrm{Berk}}} \right\rbrace
\end{equation}
Then there exists a unique isometry $W_f:L^2(\mu_f)\rightarrow \mathcal{H}_f$ of $L^2(\mu_f)$ onto $\mathcal{H}_f$ such that 
\begin{equation}
W_f(1)=f, \qquad \text{and} \qquad W_f \widehat{P}(\chi_{Z(Q)})W_f^*=S_QS_Q^*
\end{equation}
for all words $Q \in \mathcal{L}_{\mathbb{H}_{\mathrm{Berk}}^{\mathbb{Q}}(\mathbb{C}_p)}$.
\end{mcor}
\begin{proof}
To define $W_f:L^2(\mu_f)\rightarrow \mathcal{H}_f$, we set
\begin{equation*}
W_f(1)=f, \quad \mathrm{and} \quad W_f(\chi_{Z(Q)})=S_QS_Q^*f 
\end{equation*}
which then satisfies 
\begin{equation*}
\int_{\mathbb{P}^1_{\mathrm{Berk}}(\mathbb{C}_p)}|\chi_Q(x)|^2\ d\mu_f(x)=\left\langle S_QS_Q^*f, S_QS_Q^*f \right\rangle =\|W_f(\chi_Q)\|_{\mathcal{H}_f}^2
\end{equation*}
and thus extends to an isometry $W_f:L^2(\mu_f)\rightarrow \mathcal{H}_f$. Again, uniqueness follows from density of simple functions in $L^2(\mu_f)$.
\end{proof}
\subsection{Perron-Frobenius Operator}
Consider the \textit{transfer operator} $T_\sigma:L^2(\mathbb{P}^1_{\mathrm{Berk}}(\mathbb{C}_p),\mu_G)\rightarrow L^2(\mathbb{P}^1_{\mathrm{Berk}}(\mathbb{C}_p),\mu_G)$ that composes with the coding map denoted by $\sigma:\mathbb{P}^1_{\mathrm{Berk}}(\mathbb{C}_p)\rightarrow \mathbb{P}^1_{\mathrm{Berk}}(\mathbb{C}_p)$,
\begin{equation}
(T_\sigma\psi)(x)=\psi(\sigma(x)).
\end{equation}
One can associate to the operator $T_\sigma$ its adjoint $P_\sigma$ given by 
\begin{equation}
\int\psi P_\sigma(\xi)d\mu = \int T_\sigma(\psi)\xi d\mu
\end{equation}
and called the \textit{Perron-Frobenius operator}.
\begin{mth}
Let $\lbrace \sigma_{q}\rbrace$ be the semibranching system defined associated to the coding map given by $\sigma:\mathbb{P}^1_{\mathrm{Berk}}(\mathbb{C}_p)\rightarrow \mathbb{P}^1_{\mathrm{Berk}}(\mathbb{C}_p)$. Then, the Perron-Frobenius operator $P_\sigma$ is of the form
\begin{equation}
(P_\sigma\xi)(x)=\sum_{q}\chi_{Z(q)}\xi(\sigma_{q}(x)).
\end{equation}
\end{mth}
\begin{proof}
In the Hilbert space $L^2(\mathbb{P}^1_{\mathrm{Berk}}(\mathbb{C}_p),\mu_G)$, we can write
\begin{align*}
\left\langle T_\sigma(\psi),\xi\right\rangle &= \int_{\mathbb{P}^1_{\mathrm{Berk}}(\mathbb{C}_p)} \overline{\psi(\sigma(x))}\xi(x) \ d\mu(x)\\
&= \sum_{q\in \mathbb{Q}_1}  \int_{F(q)} \overline{\psi(u)}\xi(\sigma_q(u))\ d\mu(u)\\
&=\left\langle \psi, \sum_{q\in \mathbb{Q}_1} \chi_{F(q)}\xi\circ \sigma_q \right\rangle 
\end{align*}
The left-hand side is recognized as the Perron-Frobenius operator.
\end{proof}
\noindent
Since $S_q^*$ is defined as
\begin{equation*}
S_q^*\xi(x)=\chi_{F(q)}\xi\circ \sigma_q(x)
\end{equation*}
we immediately deduce the following corollary.
\begin{mcor}
Let $\lbrace \sigma_{q}\rbrace$ be the semibranching system defined associated to the coding map given by $\sigma:\mathbb{P}^1_{\mathrm{Berk}}(\mathbb{C}_p)\rightarrow \mathbb{P}^1_{\mathrm{Berk}}(\mathbb{C}_p)$. Then, the Perron-Frobenius operator $P_\sigma$ is of the form
\begin{equation}
P_\sigma=\sum_{q}S_{q}^*
\end{equation}
and $P_\sigma$ belongs to $A(\mathcal{H})$.
\end{mcor}
\section{KMS states and invariant measures}
\subsection{Unitary representations of $\mathrm{PGL}_2(\mathbb{C}_p)$}
We would like to start this section with some remarks on representations of $\mathrm{PGL}_2(\mathbb{C}_p)$ as a discrete group in the Hilbert space $L^2(\mathbb{P}^1_{\mathrm{Berk}}(\mathbb{C}_p),\mu_G)$.

\paragraph*{Gromov hyperbolic space}
Let $x,y,z$ be points of $\mathbb{P}^1_{\mathrm{Berk}}(\mathbb{C}_p)$, not all equal. We define the \textit{Gromov product} denoted $(x|y)_z$ by
\begin{equation}
(x|y)_z=\rho(w,z),
\end{equation}
where $w$ is the first point where the unique paths from $x$ to $z$ and $y$ to $z$ intersect. By convention, we set $(x|y)_z=+\infty$ if $x=y$ and $x$ is a point of type I, and we set $(x|y)_z=0$ if $x=z$ or $y=z$.\\
If $x,y,z\in \mathbb{H}_{\mathrm{Berk}}(\mathbb{C}_p)$, then one checks easily that
\begin{equation}
(x|y)_z=\frac{1}{2}\left( \rho(x,z) + \rho(y,z) -\rho(x,y) \right).
\end{equation}
This is the usual definition of the Gromov product in Gromov's theory of $\delta$-hyperbolic spaces, with $\mathbb{H}_{\mathrm{Berk}}(\mathbb{C}_p)$ being an example of a $0$-hyperbolic space.
\begin{mprop}[\cite{coornaert_mesures_1993}]
$\mathbb{H}_{\mathrm{Berk}}(\mathbb{C}_p)$ is Gromov $0$-hyperbolic.
\end{mprop}
\noindent
Let $\zeta_G$ be the Gauss point of $\mathbb{P}^1_{\mathrm{Berk}}(\mathbb{C}_p)$. Define the \textit{fundamental potential kernel relative}, written $\kappa_z(x,y)$, and the \textit{canonical distance relative to $z$}, written $[x,y]_z$, by setting
\begin{equation}
\kappa_z(x,y)=-\log_v[x,y]_z=(x|y)_\zeta - (x|z)_\zeta -(y|z)_\zeta
\end{equation} 
In particular, if one takes $z=\zeta_G$, then the fundamental kernel simplifies:
\begin{equation}
\kappa_G(x,x)=-\log_v[x,x]_G=-\log_v(\text{diam}(x)).
\end{equation}
Then, we define the following measure on $\mathbb{H}_{\mathrm{Berk}}(\mathbb{C}_p)$
\begin{equation}
d\mu (x)=p^{-\kappa_G(x,x)}d\mu_G(x)=\text{diam}(x)d\mu_G(x).
\end{equation}
Since the measure $\mu_G$ is $\mathrm{PGL}_2(\mathbb{C}_p)$-invariant, the measure $\mu$ is quasi-invariant under the action 
\begin{equation}
\sigma_g:\mathrm{PGL}_2(\mathbb{C}_p)\times \mathbb{P}^1_{\mathrm{Berk}}(\mathbb{C}_p)\rightarrow\mathbb{P}^1_{\mathrm{Berk}}(\mathbb{C}_p), \qquad \sigma_g(x)=\left( \begin{array}{cc}
z & a \\ 
0 & 1
\end{array}\right) \cdot x  
\end{equation}
The corresponding Radon-Nikodym derivatives are given by
\begin{equation}
\delta(g)(x)=\frac{d(\sigma_g)_*\mu}{d\mu}(x)=\frac{\mathrm{diam}(g\cdot x)}{\mathrm{diam}(x)}=|z|
\end{equation}
\begin{mprop}
The following maps
\begin{equation}
(U_{s}(g)f):=e^{is}\sqrt{|z|}f(g\cdot x)
\end{equation}
define a family of unitary representation of the group $\mathrm{PGL}_2(\mathbb{C}_p)$ on $L^2(\mathbb{H}_{\mathrm{Berk}}(\mathbb{C}_p),\mu)$, parametrized by $s\in \mathbb{R}$.
\end{mprop}
\begin{proof}
Clearly, $U_s(g)f$ is measurable, and we also find
\begin{align*}
\|U_s(g)f\|_2^2=\int_{\mathbb{H}_{\mathrm{Berk}}(\mathbb{C}_p)}|z| |f(g\cdot x)|^2 \ d\mu(x)&= \int_{\mathbb{H}_{\mathrm{Berk}}(\mathbb{C}_p)}|f(g\cdot x)|^2 \ \frac{d(\sigma_g)_*\mu}{d\mu} d\mu(x)\\
&=\int_{\mathbb{H}_{\mathrm{Berk}}(\mathbb{C}_p)}|f(g\cdot x)|^2 \ d(\sigma_g)_*\mu(x)\\
&=\int_{\mathbb{H}_{\mathrm{Berk}}(\mathbb{C}_p)} |f(x)|^2 \ d\mu(x)=\|f\|^2
\end{align*}
Thus, $U_s(g)f$ defines an isometry of $L^2(\mathbb{H}_{\mathrm{Berk}}(\mathbb{C}_p),\mu)$. We also observe that, for $g,h\in \mathrm{PGL}_2(\mathbb{C}_p)$, we can write
\begin{equation}
U_s(gh)f= e^{2is}\sqrt{|z_g|}\sqrt{|z_h|}f(gh\cdot x)=U_s(g)\left( U_s(h)f\right) 
\end{equation}
In particular, we see that each isometry $U_s(g)$ is surjective with $U_s(g^{-1})=U_s(g)^{-1}$.
\end{proof}
\noindent
Let us recall the definition of a covariant representation of a triple $C^*$-dynamical system.
\begin{mdef}
Let $\alpha:G\rightarrow \mathrm{Aut}(A)$ be an action of a locally compact group $G$ on a $C^*$-algebra $A$. A \textit{covariant representation} of $(G,A,\alpha)$ on a Hilbert space $H$ is a pair $(v,\pi)$ consisting of a unitary representation $v:G\rightarrow U(H)$ (the unitary group of $H$) and a representation $\pi:A\rightarrow B(H)$, satisfying  the \textit{covariant condition}
\begin{equation*}
v(g)\pi(a)v(g)^*=\pi(\alpha_g(a))
\end{equation*}
for all $g\in G$ and $a\in A$. It is called \textit{nondegenerate} if $\pi$ is nondegenerate.
\end{mdef}
\noindent
In the present context, we consider the $C^*$-algebra $C_0(\mathbb{H}_{\mathrm{Berk}}(\mathbb{C}_p))$ of continuous functions on the hyperbolic Berkovich line. We consider the representation by left-multiplication :
\begin{equation}
\pi: C_0(\mathbb{H}_{\mathrm{Berk}}(\mathbb{C}_p))\rightarrow B(L^2(\mathbb{H}_{\mathrm{Berk}}(\mathbb{C}_p),\mu)), \qquad \pi(f)\psi(x)=f(x)\psi(x)
\end{equation}
and we will denote the Hilbert space $L^2(\mathbb{H}_{\mathrm{Berk}}(\mathbb{C}_p)$ by $\mathcal{H}$.
\begin{mth}
For $s\in \mathbb{R}$, the pair $(U_s,\pi)$ is a \textit{covariant representation} of the following triple \\ $(\mathrm{PGL}_2(\mathbb{C}_p),C_0(\mathbb{H}_{\mathrm{Berk}}(\mathbb{C}_p)),\sigma)$.
\end{mth}
\begin{proof}
We verify that the pair $(U_s,\pi)$ satisfies the covariant condition:
\begin{align*}
(U_s(g)\pi(a) U_s(g)^*)f(x)&=(U_s(g)\pi(a))(e^{is}|z|^{-\frac{1}{2}}f(g^{-1}\cdot x))\\
&=U_s(g)\left(e^{is}|z|^{-\frac{1}{2}}a(x)f(g^{-1}\cdot x) \right) =a(g\cdot x)f(x) 
\end{align*}
for $a\in C_0(\mathbb{H}_{\mathrm{Berk}}(\mathbb{C}_p))$ and $g\in \mathrm{PGL}_2(\mathbb{C}_p)$.
\end{proof}
\noindent
Let us denote by $\nu$ the Haar measure on $\mathrm{PGL}_2(\mathbb{C}_p)$. We let $C_c(\mathrm{PGL}_2(\mathbb{C}_p),C_0(\mathbb{H}_{\mathrm{Berk}}(\mathbb{C}_p)),\sigma)$ be the $*$-algebra of compactly supported continuous functions $f:\mathrm{PGL}_2(\mathbb{C}_p)\rightarrow C_0(\mathbb{H}_{\mathrm{Berk}}(\mathbb{C}_p))$, with pointwise addition and scalar multiplication. We define a norm $\|\cdot\|_1$ on the algebra $C_c(\mathrm{PGL}_2(\mathbb{C}_p),C_0(\mathbb{H}_{\mathrm{Berk}}(\mathbb{C}_p)),\sigma)$ by $\|f\|_1=\int_{\mathrm{PGL}_2(\mathbb{C}_p)}\|f(g)\|d\nu(g)$. Then, the space $L^1(\mathrm{PGL}_2(\mathbb{C}_p),C_0(\mathbb{H}_{\mathrm{Berk}}(\mathbb{C}_p)),\sigma)$ is the Banach $*$-algebra obtained by completion of the algebra $C_c(\mathrm{PGL}_2(\mathbb{C}_p),C_0(\mathbb{H}_{\mathrm{Berk}}(\mathbb{C}_p)),\sigma)$ with respect to $\|\cdot\|_1$.\\

\noindent
We can define a representation $\alpha:C_c(\mathrm{PGL}_2(\mathbb{C}_p),C_0(\mathbb{H}_{\mathrm{Berk}}(\mathbb{C}_p)),\sigma)\rightarrow B(H)$ given by
\begin{equation}
\alpha_s(f)\psi(x) = \int_{\mathrm{PGL}_2(\mathbb{C}_p)}\pi(f_g)U_s(g)\psi d\nu(g)=\int_{\mathrm{PGL}_2(\mathbb{C}_p)}f_g(x)\psi(g\cdot x)e^{is}\sqrt{|z|}d\nu(g)
\end{equation}
and define the \textit{crossed product} $C^*$-algebra $C_0(\mathbb{H}_{\mathrm{Berk}}(\mathbb{C}_p))\rtimes_\sigma \mathrm{PGL}_2(\mathbb{C}_p)$ to be the norm closure of the image $\alpha_s(L^1(\mathrm{PGL}_2(\mathbb{C}_p),C_0(\mathbb{H}_{\mathrm{Berk}}(\mathbb{C}_p)),\sigma))$.
\begin{mprop}
The pair $(\alpha_s,H)$ defines a representation of the crossed-product $C^*$-algebra $C_0(\mathbb{H}_{\mathrm{Berk}}(\mathbb{C}_p))\rtimes_\sigma \mathrm{PGL}_2(\mathbb{C}_p)$.
\end{mprop}
\subsection{K-cycle and dynamical system}
In the previous section, we have constructed a $C^*$-algebra on the hyperbolic projective line $\mathbb{H}_{\mathrm{Berk}}(\mathbb{C}_p)$. If we see $\mathbb{P}^1_{\mathrm{Berk}}(\mathbb{C}_p)$ as a tree, then we studied crossed product $C^*$-algebra on the interior of the tree. In this last section, we show that crossed product $C^*$-algebra can also be used to study the boundary $\partial\mathbb{P}^1_{\mathrm{Berk}}(\mathbb{C}_p)$ which is identified with $\mathbb{P}^1(\mathbb{C}_p)$.\\

\noindent
Consider $\mathrm{\Gamma}$ be a finitely generated subgroup of $\mathrm{PGL}_2(\mathbb{C}_p)$. Moreover, let $\mathrm{\Gamma}$ be torsion-free, in which case it is called a (p-adic) \textit{Schottky group}. In fact, a Schottky group is a free group \cite{maskit_characterization_1967}.
\begin{mdef}[Limit set, \cite{hersonsky_groups_1997}]
Let $(X,d)$ be a metric space and $\mathrm{\Gamma}\subset \text{Iso}(X)$ be a discrete subgroup. The \textit{limit set} of $\mathrm{\Gamma}$ is the subset $\mathrm{\Lambda}=\mathrm{\Lambda}(\mathrm{\Gamma})$ of $\partial X$ of points which are accumulation points of orbits in $X$. That is 
\begin{equation*}
\mathrm{\Lambda} := \Bigl\{ y\in \partial X \ | \ y=\lim\lbrace \gamma_m(x)\rbrace \ \text{for some $x\in X$ and $\lbrace \gamma_m\rbrace$ a sequence in Iso$(X)$}\Bigr\} 
\end{equation*}
\end{mdef}
\noindent
Since $\mathrm{PGL}_2(\mathbb{C}_p)$ acts transitively on type II points in $\mathbb{P}^1_{\mathrm{Berk}}(\mathbb{C}_p)$, then the limit set of $\mathrm{\Gamma}$ is $\mathbb{P}^1(\mathbb{C}_p)$. Starting from the commutative algebra of Lipschitz-continuous functions $C(\mathbb{P}^1(\mathbb{C}_p))$ over $\mathbb{P}^1(\mathbb{C}_p)$, we will consider the \textit{reduced cross-product} for which we briefly recall the construction. Consider the integrated representation $\alpha$ define in the previous section:
\begin{equation*}
\alpha:L^1(\mathrm{PGL}_2(\mathbb{C}_p),C(\mathbb{P}^1),\sigma)\rightarrow B(\mathcal{H}),\qquad \alpha(f)\psi(\xi)=\sum_{\gamma\in \mathrm{\Gamma}}f_\gamma(\xi) (U_\gamma \psi)(\xi)
\end{equation*}
\noindent
The \textit{reduced cross-product $C^*$-algebra} $C^*_r(\mathrm{\Gamma},\mathbb{P}^1(\mathbb{C}_p))$ is the completion of the Banach algebra $L^1(\mathrm{PGL}_2(\mathbb{C}_p),C(\mathbb{P}^1(\mathbb{C}_p)))$ with respect to the operator norm
\begin{equation}
f\mapsto \|\alpha(f)\|_{B(\mathcal{H})}.
\end{equation} 
Implicitly, the definition of $C^*_r(\mathrm{\Gamma},\mathbb{P}^1(\mathbb{C}_p))$ is a representation of the Banach algebra \par
\noindent
$L^1(\mathrm{PGL}_2(\mathbb{C}_p),C(\mathbb{P}^1(\mathbb{C}_p)))$ and therefore of the full cross-product $C(\mathbb{P}^1(\mathbb{C}_p))\rtimes_\sigma \mathrm{PGL}_2(\mathbb{C}_p)$. Thus, there is a surjective *-homomorphism between the crossed product algebras $C(\mathbb{P}^1(\mathbb{C}_p))\rtimes_\sigma \mathrm{PGL}_2(\mathbb{C}_p)\rightarrow C^*_r(\mathrm{\Gamma},\mathbb{P}^1(\mathbb{C}_p))$.
\paragraph{Time evolution} 
One can define the time evolution as the one-parameter family acting on the reduced $C^*$-algebra $C^*_r(\mathrm{\Gamma},\mathbb{P}^1(\mathbb{C}_p))$ given by
\begin{equation}
\alpha_t\left( \sum_\gamma f_\gamma(\xi)U_\gamma\right)=\sum_\gamma p^{itB_\xi(\zeta,\gamma\cdot\zeta)}f_\gamma(\xi)U_\gamma
\end{equation}
where $B(\zeta,\gamma \zeta,\xi)$ is the \textit{Busemann function} based at the Gauss point $\zeta$ and, defined by:
\begin{equation}
B(\zeta,\gamma\zeta,\xi)=\lim_{\substack{x\rightarrow \xi\\ x\in \mathbb{H}^1_{\mathrm{Berk}}}} (\zeta,x)_{\gamma\zeta} - (\gamma\zeta, x)_{\gamma} = \lim_{\substack{x\rightarrow \xi\\ x\in \mathbb{H}^1_{\mathrm{Berk}}}} \log_v\left( \frac{\mathrm{diam}(\gamma^{-1}x)}{ \mathrm{diam}(x)}\right)
\end{equation}
for all $\xi \in \mathbb{P}^1(\mathbb{C}_p)$, using the Gromov product.
\begin{mlem}
\label{lemunit}
$\alpha$ is a strongly-continuous one parameter group of $*$-automorphisms of the cross-product algebra $C(\mathbb{P}^1_{\mathrm{Berk}}(\mathbb{C}_p)
)\rtimes \mathrm{\Gamma}$.
\end{mlem}
\begin{proof}
Let $\gamma' \mapsto (\xi \mapsto h(\xi)U_{\gamma'})\in \mathcal{H}$, then define the unitary operator $v(t,\xi)$ by
\begin{equation}
v(t,\xi)h(\xi)U_{\gamma'} = p^{it B_\xi(\zeta,\gamma'\cdot \zeta)}h(\xi)U_{\gamma'}
\end{equation}
Then, consider the operators  $v(t,\xi)f(\xi)U_\gamma v(t,\xi)^*$ which acts on the Hilbert space $\mathcal{H}$ as follows 
\begin{align*}
\left[ v(t,\xi)f(\xi)U_\gamma v(t,\xi)^*\right] h(\xi)U_{\gamma'}&=\left[ v(t,\xi)f(\xi)\right]p^{-it B_{\gamma^{-1}\xi}(\zeta,\gamma'\cdot \zeta)}h(\gamma^{-1}\cdot \xi)U_{\gamma\gamma'}\\
&=p^{it(B_\xi(\zeta,\gamma\gamma'\cdot \zeta)-B_{\gamma^{-1}\xi}(\zeta,\gamma'\cdot \zeta))}f(\xi)h(\gamma^{-1}\cdot \xi)U_{\gamma\gamma'}
\end{align*}
From the definition of the Busemann function, we can show that 
\begin{equation*}
-B_{\gamma^{-1}\xi}(\zeta,\gamma'\cdot \zeta)=B_{\gamma^{-1}\xi}(\gamma'\cdot \zeta, \zeta), \quad B_{\gamma^{-1}\xi}(\zeta,\gamma'\cdot \zeta)=B_{\xi}(\gamma\zeta,\gamma\gamma'\cdot \zeta)
\end{equation*}
and $B_\xi(\zeta,\gamma\gamma'\cdot \zeta) + B_{\xi}(\gamma\gamma'\cdot \zeta, \gamma\zeta)= B_\xi(\zeta,\gamma\zeta )$, from which we deduce that 
\begin{equation*}
\left[ v(t,\xi)f(\xi)U_\gamma v(t,\xi)^*\right] h(\xi)U_{\gamma'}=p^{itB_\xi(\zeta, \gamma\cdot \zeta)}f(\xi)h(\gamma^{-1}\cdot \xi)U_{\gamma\gamma'}=\alpha_t\left( f(\xi)U_\gamma\right) h(\xi)U_{\gamma'}
\end{equation*}
Thus, the map $\alpha_t$ coincides on the reduced $C^*$-algebra $C^*_r(\mathrm{\Gamma},\mathbb{P}^1(\mathbb{C}_p))$ with the conjugate by the unitary operator $v(t,\xi)$. Hence, $t\mapsto\alpha_t$ is a strongly continuous $*$-automorphism which clearly satisfies for any time $t,t'\in \mathbb{R}$, $\alpha_t\circ \alpha_{t'}=\alpha_{t+t'}$.
\end{proof}
\noindent
\paragraph*{The Patterson-Sullivan measure on $\mathbb{P}^1_{\mathrm{Berk}}(\mathbb{C}_p)$}
Let us recall the definition of the Patterson-Sullivan measure. In addition to the Busemann function $B(\zeta,\gamma\zeta,\xi)$, consider the function
\begin{equation}
j_\gamma(\xi)=p^{B(\zeta,\gamma\zeta,\xi)}=\lim_{\substack{x\rightarrow \xi\\ x\in \mathbb{H}^1_{\mathrm{Berk}}}}  \frac{\mathrm{diam}(\gamma^{-1}x)}{ \mathrm{diam}(x)}
\end{equation}
\begin{mdef}
Let $D\geq 0$ and $\mu$ a finite (and non-trivial) measure on $\mathbb{P}^1(\mathbb{C}_p)$. The measure $\mu$ is called $\mathrm{\Gamma}$-quasiconformal of dimension $D$ if the pullback measures $\gamma^*\mu$ ($\gamma\in \mathrm{\Gamma}$) are pairwise absolutely continuous and there exists $C\geq 1$ such that
\begin{equation}
C^{-1}j_\gamma(\xi)^D\leq \frac{d(\gamma^*\mu)}{d\mu}\leq Cj_\gamma^D, \qquad \mu-a.e.
\end{equation}
\end{mdef}
\noindent
Let us now recall the definition of the \textit{critical exponent} of $\Gamma$ as the number
\begin{equation}
\delta(\mathrm{\Gamma})=\limsup_{R\rightarrow \infty}\frac{1}{R}\log \left\lbrace \gamma \in \mathrm{\Gamma} : \rho(\zeta,\gamma\zeta)\leq R\right\rbrace 
\end{equation}
It is equivalently defined as the abscissa of convergence of the Poincaré series
\begin{equation}
\wp_{\mathrm{\Gamma}}(s)=\sum_{\gamma\in \mathrm{\Gamma}}e^{-s\mathrm{diam}(\gamma\zeta)}
\end{equation}
\begin{mprop}[Theorem 5.4 in \cite{coornaert_mesures_1993}]
If $\delta(\mathrm{\Gamma})<\infty$, then there exists a (unique) measure $\mathrm{\Gamma}$-quasiconformal of dimension $\delta(\mathrm{\Gamma})$ with support coinciding with $\mathbb{P}^1(\mathbb{C}_p)$. This measure is called the Patterson-Sullivan measure and denoted by $\mu_{PS,\zeta}$.
\end{mprop}
\noindent
We can then define the positive functional $\tau:C^*_r(\mathrm{\Gamma},\mathbb{P}^1(\mathbb{C}_p))\rightarrow \mathbb{C}$ by 
\begin{equation}
\tau(f)=\int_{\mathbb{P}^1(\mathbb{C}_p)}f_e(\xi)d\mu_{PS,\zeta}(\xi)
\end{equation}
\paragraph*{$\mathrm{KMS}_\beta$-states.}
Consider a $C^*$-dynamical i.e. a pair $(A,\alpha)$ with $A$ a $C^*$-algebra and $\alpha:\mathbb{R}\rightarrow \mathrm{Aut}(A)$ a time evolution.
\begin{mdef}[KMS state]
A state $\varphi$ on $A$ satisfies the Kubo-Martin-Schwinger (KMS) condition with respect to $\alpha$ at inverse temperature $\beta\in [0,\infty)$ ($\varphi$ is a $\alpha$-$\mathrm{KMS}_\beta$ state), if 
\begin{equation}
\varphi(ab)=\varphi(b\alpha_{i\beta}(a))
\end{equation}
for all element $a,b$ in a dense subalgebra $A^{an}\subset A$ of \textit{analytic elements}.
\end{mdef}
\noindent
Then, we can show that the $C^*$-dynamical $(C^*_r(\mathrm{\Gamma},\mathbb{P}^1(\mathbb{C}_p)), \alpha_t)$ admits a unique KMS-state that is exactly given by the Patterson-Sullivan measure on the boundary of the Berkovich projective line.
\begin{mprop}[Proposition 5.11 \cite{lott_limit_2005}]
\label{thm3}
The reduced cross-product $C^*$-algebra $C^*_r(\mathrm{\Gamma},\mathbb{P}^1(\mathbb{C}_p))$ has a unique $\mathrm{KMS}_{\beta}$ states with inverse temperature $\beta=\delta(\mathrm{\Gamma})$ obtained by integration with the Patterson-Sullivan measure:
\begin{equation}
\varphi_{\beta,\zeta}\left(  \sum_\gamma f_\gamma(\xi)U_\gamma\right) =\int_{\mathbb{P}^1(\mathbb{C}_p)}f_e(\xi)d\mu_{PS,\zeta}(\xi).
\end{equation}
\end{mprop}
\paragraph*{Hamiltonian operator}
If we denote the family Radon-Nikodym derivatives as, 
\begin{equation}
\delta(\gamma)=\frac{d\gamma^*\mu_{PS}}{d\mu_{PS}}
\end{equation}
then as previously the family of operators acting on $L^2(\mathbb{P}^1,\mu_{PS})$:
\begin{equation}
U(\gamma)f(\xi)=\sqrt{\delta(\gamma)}f(\gamma\cdot \xi)
\end{equation}
defines a unitary representation of the Schottky group $\mathrm{\Gamma}$. Again we denote by $\pi$ the representation by left multiplication of $C(\mathbb{P}^1(\mathbb{C}_p))$ on $L^2(\mathbb{P}^1(\mathbb{C}_p),\mu_{PS})$.
\begin{mth}
The pair $(\pi,U)$ of maps defines a regular covariant representation of the crossed product $C^*$-algebra $C(\mathbb{P}^1(\mathbb{C}_p))\rtimes \mathrm{\Gamma}$. In this representation, the time evolution is implemented by the Hamiltonian 
\begin{equation}
H f(\xi)U(\gamma) = B(\zeta,\gamma\zeta,\xi)f(\xi)U(\gamma)
\qquad
\text{such that}
\qquad 
\alpha_t(a)=e^{itH}ae^{-itH}
\end{equation}
\end{mth}
\begin{proof}
From Lemma \ref{lemunit}, we have shown the time evolution $\alpha_t$ is given as the conjugate of a unitary operator $v(t,\xi)$ such that 
\begin{equation*}
\alpha_t(a)(\xi)=v(t,\xi)a(\xi)v(t,\xi)^*
\end{equation*}
Then, the statement follows from the fact that there exists a selfadjoint operator $H$ acting on the Hilbert space $\mathcal{H}$ such that $v(t,\xi)=e^{itH}$.
\end{proof}
\noindent
Following \cite{greenfield_twisted_2014}, we can then introduce on the boundary $\partial\mathbb{P}^1_{\mathrm{Berk}}(\mathbb{C_p})$ a (twisted) spectral triple given by  $(C^*_r(\mathrm{\Gamma},\mathbb{P}^1(\mathbb{C}_p)), \mathcal{H}, D)$ where the Dirac operator $D$ is such that $|D|=H$.
\addcontentsline{toc}{section}{References}
\bibliographystyle{amsplain}
\bibliography{References.bib}

\providecommand{\bysame}{\leavevmode\hbox to3em{\hrulefill}\thinspace}
\providecommand{\MR}{\relax\ifhmode\unskip\space\fi MR }
\providecommand{\MRhref}[2]{%
  \href{http://www.ams.org/mathscinet-getitem?mr=#1}{#2}
}
\providecommand{\href}[2]{#2}
\begin{thebibliography}{10}

\bibitem{baker_-adic_2008}
Matthew Baker, Brian Conrad, Samit Dasgupta, Kiran Kedlaya, and Jeremy
  Teitelbaum, \emph{p-adic {Geometry}}, University {Lecture} {Series}, vol.~45,
  American Mathematical Society, Providence, Rhode Island, August 2008 (en).

\bibitem{baker_potential_nodate}
Matthew Baker and Robert~S. Rumely, \emph{Potential theory and dynamics on the
  {Berkovich} projective line}, Mathematical surveys and monographs, no. v.
  159, American Mathematical Society, Providence, R.I, 2010.

\bibitem{banic_wazewskis_2013}
Iztok Banič, Matevž Črepnjak, Matej Merhar, Uroš Milutinović, and Tina
  Sovič, \emph{Ważewski's universal dendrite as an inverse limit with one
  set-valued bonding function}, Glasnik Matematicki \textbf{48} (2013), no.~1,
  137--165 (en).

\bibitem{boava_c-algebras_2023}
Giuliano Boava, Gilles~G. de~Castro, Daniel Gonçalves, and Daniel~W. van Wyk,
  \emph{C*-{Algebras} of one-sided subshifts over arbitrary alphabets},
  December 2023, arXiv:2312.17644 [math].

\bibitem{boava_algebras_2023}
Giuliano Boava, Gilles~G. De~Castro, Daniel Gonçalves, and Daniel~W. Van~Wyk,
  \emph{Algebras of one-sided subshifts over arbitrary alphabets}, Revista
  Matemática Iberoamericana \textbf{40} (2024), no.~3, 1045--1088 (en).

\bibitem{bratteli_iterated_1996}
Ola Bratteli and Palle E.~T. Jørgensen, \emph{Iterated function systems and
  permutation representations of the {Cuntz} algebra}, American {Mathematical}
  {Society}, no. volume 139, number 663, American Mathematical Society, 1999
  (eng).

\bibitem{charatonik_self-homeomorphic_1994}
Wodzimierz~J. Charatonik and Anne Dilks, \emph{On self-homeomorphic spaces},
  Topology and its Applications \textbf{55} (1994), no.~3, 215--238 (en).

\bibitem{connes_non-commutative_1985}
A.~Connes, \emph{Noncommutative differential geometry}, Publications
  mathématiques de l'IHÉS \textbf{62} (1985), no.~1, 41--144 (en).

\bibitem{connes_noncommutative_1994}
\bysame, \emph{Noncommutative geometry}, Academic Press, San Diego, 1994 (eng).

\bibitem{connes_noncommutative_2007}
A.~Connes and M.~Marcolli, \emph{Noncommutative {Geometry}, {Quantum} {Fields}
  and {Motives}}, Colloquium {Publications}, vol.~55, American Mathematical
  Society, Providence, Rhode Island, December 2007 (en).

\bibitem{connes_scaling_2015}
Alain Connes and Caterina Consani, \emph{The scaling site}, Comptes Rendus.
  Mathématique \textbf{354} (2015), no.~1, 1--6 (en).

\bibitem{connes_geometry_2016}
\bysame, \emph{Geometry of the arithmetic site}, Advances in Mathematics
  \textbf{291} (2016), 274--329 (en).

\bibitem{connes_physics_nodate}
Alain Connes and Matilde Marcolli, \emph{From {Physics} to {Number} {Theory}
  via {Noncommutative} {Geometry}} (en).

\bibitem{consani_spectral_2003}
Caterina Consani and Matilde Marcolli, \emph{Spectral triples from {Mumford}
  curves}, International Mathematics Research Notices \textbf{2003} (2003),
  no.~36, 1945 (en).

\bibitem{consani_non-commutative_2003}
\bysame, \emph{Noncommutative geometry, dynamics, and $\infty$-adic {Arakelov}
  geometry}, Selecta Mathematica \textbf{10} (2004), no.~2, 167--251 (en).

\bibitem{coornaert_mesures_1993}
Michel Coornaert, \emph{Mesures de {Patterson}-{Sullivan} sur le bord d’un
  espace hyperbolique au sens de {Gromov}}, Pacific Journal of Mathematics
  \textbf{159} (1993), no.~2, 241--270 (fr).

\bibitem{CORNELISSEN2013110}
Gunther Cornelissen and Matilde Marcolli, \emph{Graph reconstruction and
  quantum statistical mechanics}, Journal of Geometry and Physics \textbf{72}
  (2013), 110--117, Noncommutative algebraic geometry and its applications to
  physics.

\bibitem{dutkay_iterated_2006}
Dorin Dutkay and Palle Jorgensen, \emph{Iterated function systems, {Ruelle}
  operators, and invariant projective measures}, Mathematics of Computation
  \textbf{75} (2006), no.~256, 1931--1970 (en).

\bibitem{dutkay_atomic_2015}
Dorin~Ervin Dutkay, John Haussermann, and Palle~E.T. Jorgensen, \emph{Atomic
  representations of {Cuntz} algebras}, Journal of Mathematical Analysis and
  Applications \textbf{421} (2015), no.~1, 215--243 (en).

\bibitem{exel_cuntz-krieger_1999}
R.~Exel and M.~Laca, \emph{Cuntz-{Krieger} algebras for infinite matrices},
  Journal für die reine und angewandte Mathematik (Crelles Journal)
  \textbf{1999} (1999), no.~512, 119--172.

\bibitem{floricel_inductive_2017}
Remus Floricel and Asghar Ghorbanpour, \emph{On inductive limit spectral
  triples}, Proceedings of the American Mathematical Society (2017).

\bibitem{fremlin_measure_2006}
David~H. Fremlin, \emph{Measure theory. 4,2: {Topological} measure spaces}, 2.
  print ed., Torres Fremlin, Colchester, 2006 (en).

\bibitem{greenfield_twisted_2014}
M.~Greenfield, M.~Marcolli, and K.~Teh, \emph{Twisted spectral triples and
  quantum statistical mechanical systems}, P-Adic Numbers, Ultrametric
  Analysis, and Applications \textbf{6} (2014), no.~2, 81--104 (en).

\bibitem{gubser_edge_2017}
Steven~S. Gubser, Matthew Heydeman, Christian Jepsen, Matilde Marcolli, Sarthak
  Parikh, Ingmar Saberi, Bogdan Stoica, and Brian Trundy, \emph{Edge length
  dynamics on graphs with applications to p-adic {AdS}/{CFT}}, Journal of High
  Energy Physics \textbf{2017} (2017), no.~6, 157 (en).

\bibitem{hersonsky_groups_1997}
Sa'Ar Hersonsky and John Hubbard, \emph{Groups of automorphisms of trees and
  their limit sets}, Ergodic Theory and Dynamical Systems \textbf{17} (1997),
  no.~4, 869--884 (en).

\bibitem{heydeman_tensor_2017}
Matthew Heydeman, Matilde Marcolli, Ingmar Saberi, and Bogdan Stoica,
  \emph{Tensor networks, \$p\$-adic fields, and algebraic curves: arithmetic
  and the {AdS}\$\_3\$/{CFT}\$\_2\$ correspondence}, Adv. Theor. Math. Phys.
  \textbf{22} (2018), no.~1 (en).

\bibitem{jorgensen_states_2011}
P.~E.~T. Jorgensen and A.~M. Paolucci, \emph{{States} {on} {the} {Cuntz}
  {Algebras} {and} \textit{p} -{adic} {random} {walks}}, Journal of the
  Australian Mathematical Society \textbf{90} (2011), no.~2, 197--211 (en).

\bibitem{jorgensen_use_2005}
Palle E.~T. Jorgensen, \emph{Use of operator algebras in the analysis of
  measures from wavelets and iterated function systems}, arXiv: Operator
  Algebras (2005).

\bibitem{kajiwara_c-algebras_2003}
Tsuyoshi Kajiwara and Yasuo Watatani, \emph{${C}^*$-{Algebras} {Associated}
  with {Complex} {Dynamical} {Systems} and {Backward} {Orbit} {Structure}},
  January 2014, pp.~243--254.

\bibitem{lott_limit_2005}
John Lott, \emph{Limit {Sets} as {Examples} in {Noncommutative} {Geometry}},
  K-Theory \textbf{34} (2005), no.~4, 283--326 (en).

\bibitem{marcolli_cuntzkrieger_2011}
Matilde Marcolli and Anna~Maria Paolucci, \emph{Cuntz–{Krieger} {Algebras}
  and {Wavelets} on {Fractals}}, Complex Analysis and Operator Theory
  \textbf{5} (2011), no.~1, 41--81 (en).

\bibitem{maskit_characterization_1967}
Bernard Maskit, \emph{A characterization of {Schottky} groups}, Journal
  d'Analyse Mathématique \textbf{19} (1967), no.~1, 227--230 (en).

\bibitem{Poineau2021}
J{\'e}r{\^o}me Poineau and Daniele Turchetti, \emph{Berkovich curves and
  schottky uniformization i: The berkovich affine line}, pp.~179--223, Springer
  International Publishing, Cham, 2021.

\bibitem{ren_linear_nodate}
Qiuyu Ren, \emph{Linear {Fractional} {Transformations} on the {Berkovich}
  {Projective} {Line}} (en).

\bibitem{stoll_p-adic_nodate}
Michael Stoll, \emph{p-adic {Analysis} in {Arithmetic} {Geometry}} (en).

\end{thebibliography}
\end{document}